# Derivation of a 2D PCCU-AENO method for nonconservative problems. Theory, Method and theoretical arguments

NGATCHA NDENGNA ARNO ROLAND

E3M laboratory, National Higher Polytechnic school of Douala, University of Douala, Cameroon

Corresponding author: arnongatcha@gmail.com

**Abstract**

In this paper, we introduce a methodology to design genuinely two-dimensional (2D) second-order path-conservative central-upwind (PCCU) schemes. The scheme studies dam-break with high sediment concentration over abrupt moving topography quickly spatially variable even in the presence of resonance. This study is possible via a 2D sediment transport model (including arbitrarily sloping sediment beds and associated energy and entropy) in new generalized Shallow Water equations derived with associated energy and entropy in this work. We establish an existence theorem of global weak solutions. We show the convergence of a sequence of solutions of the proposed model. The second-order accuracy of the PCCU scheme is achieved using a new extension AENO (Averaging Essentially Non-Oscillatory) reconstruction developed in the 2D version of this work. We prove by rigorous demonstrations that the derived 2D scheme on structured meshes is well-balanced and positivity-preserving. Several tests are made to show the ability and superb performance of the proposed numerical modeling. The results obtained are compared with those existing in the literature and with experimental data. The current modeling improves some recent results in sediment transport and shows a good ability to simulate sediment transport in large-range environments.

## I. Introduction

This work proposes a new second order finite volume method to solve a new averaged hyperbolic sediment transport model (STM) that includes arbitrarily sloping sediment beds for application in coastal or estuarine environments.

### *Sediment transport models*

Sediment transport models can be based on Saint-Venant equations, homogeneous Shallow Water equations, or nonhomogeneous Shallow Water equations. Classical sediment transport models based on homogeneous shallow water equations and the Exner model use empiric or heuristic bedload sediment flux formulas and do not able to capture internal topography waves. It's possible during the evolution of topography to observe resonance phenomena. We also can observe the situation where the flow is near the resonance. The resonance phenomenon appears when the free internal wavelength that satisfies the unforced equations coincides with the wavelength of topography forcing. The presence of the resonance permits to obtain a hypersurface on which all the characteristic fields linearly degenerate.

This phenomenon is completely ignored in several nonhomogeneous or homogeneous Shallow Water based models recently developed and which state that *the sediment velocity is equal to fluid velocity* and the topography moves with the fluid velocity [1], [2], [3] [4]. In subcritical and supercritical flow conditions these statements are not applicable. All these shortcomings make the Exner-based models a very limited model to describe the morphodynamics with accuracy. A new bedload sediment transport model that captures bed waves and accounts for the phase lag effect is proposed.

Note that when the bed moves, the classical Exner equation is not enough to properly describe the morphodynamic evolution of the channel (regular or irregular). To control the local velocity of sediment and more generally the characteristic velocity of the advection of the bed sediment form, a non-heuristic formula is used. This term corresponds to the impulse of the entrained mass that must instantly assume the characteristic velocity of the moving at the bed interface. Here, the alternative formulation of the bed evolution equation proposed extends the classical Exner model and applies to a wide range of environmental contexts.

*Hyperbolicity and mathematical analysis of sediment transport models*

From a mathematical point of view, the two-dimensional averaged sediment transport models developed in the literature admit two major difficulties related to the hyperbolicity study and mathematical analysis. It is difficult to show the existence of entropy solutions and the regularity and uniqueness of weak solution when they exist. For some sediment transport available in the literature this important part is often neglected. A rigorous mathematical analysis of a sediment transport model is performed by Birnir and Rowlett [5]. In this work, a brief mathematical analysis of the model is presented. We expose some important results.

Hyperbolicity can fail due to the morphodynamic equation used that can require complex sediment transport flux formula. To address the hyperbolicity of ST models, there are some alternatives used in the literature. Some authors use the Lagrange theorem [4] or Gerschgorin theorem [6] to find the eigenvalues of the STM when the Exner equation (integrating Grass or other complicated formulas) is used. The finding of eigenvalues of the sediment transport system of equations can depend on the choice of empiric sediment flux formula used [7] or the bedload model used [8].

It's also possible to use a splitting flux technique (which can fail in some situations) as in [9]. With splitting flux, the system becomes hyperbolic or weakly hyperbolic and the eigenstructure can be easily found without the use of Lagrange or Gerschgorin theorems. For some ST models, when the eigenvalues cannot be explicitly calculated, we use a decoupled approach for solving the problem. Such a technique is used by [10]. Due to strong and quick interactions between the flow and the moving topography, the coupled approach is often used in the literature. This technique is most appropriate than the decoupled approach that can reduce the number of total waves involved in the physic of the model. We show in [2] and [4] that the decoupled approach may fail, producing unphysical instabilities. The question of hyperbolicity study remains open for several sediment transport models when the bed evolution equations are complex. The difficulty to have a genuinely 2D hyperbolic without any ad-hoc assumptions remains for several scientists. A simple non-heuristic bedload equation is proposed here according to the kinematic equation of the bed interface to address this issue.

*Numerical schemes and limitations*

The proposed nonconservative model is addressed by a finite volume method (FVM) with special reconstruction procedures. FVM is an important building block of numerical methods for hyperbolic systems. Numerically, the formal consistency with a particular definition of weak solutions does not imply that the limits of the numerical approximations are weak solutions this major difficulty can appear with the presence of large shocks which do not satisfy the jump condition for the definition of weak solutions. Some numerical methods have been developed to solve sediment transport problems (see for instance [11], [1], [10] ). Flux-limiter scheme based on the Lax–Wendroff method coupled with a non-homogeneous Riemann solver and a flux limiter function developed by [12], needs an explicit knowledge of the eigenstructure of the system. This makes the Flux-limiter scheme to be an expensive scheme from the computational point of view and less expensive and more accurate schemes are still desired. A well-balanced positive HLLC-based scheme has been developed by Castro et al., [13]. This scheme requires an increasing number of intermediate waves and can become computationally expensive and even complex to study 2D sediment transport problems. Central-Upwind (CU) scheme or upwind numerical method requires the knowledge of the eigenstructure of a problem and often suffers from instability and robustness problems when the bed load integrates complex empiric formula. Roe-type methods based on a special linearization of the nonlinear system of governing PDE often account for all the intermediate waves and require also an explicit knowledge of the eigenstructure of the system. This makes the Roe-based method computationally expensive. Riemann HLL (Harten-Levy-Lax) solver [14] is often solicited for use in solving ST problems [15]. HLL solver is an incomplete Riemann solver and accounts only for the fastest and slowest speeds of propagation. One major drawback of the Riemann HLL solver is the increase of numerical diffusion (or dissipation). Its variants as HLLC of Toro, Spruce, and Speares [16] and HLLEM [17]  have less dissipation but require other spectral information. The use of the HLL Riemann solver to evaluate the flux is possible but can be difficult when the number of intermediate wave increase. The use of an HLLC solver can require a resolution of complex nonlinear problems via the Newton method and integrate some empiric considerations or choice of functions (see **[18]**).

More general path-conservative incomplete Riemann schemes or its extension can also be used for sediment transport. Amounts these schemes we cite the PVM (polynomial viscosity matrix) and RVM (Rational viscosity matrix) of Castro et al., **[19]**, **[20]** or both PVM and RVM solvers and their variants (see for instance **[18]**). All these schemes require a choice of function to control the numerical diffusion and some other empiric considerations whose designing is not easy. Based on a formalism of path-conservative **[21]**, PCCU has been designed to improve some classical nonconservative schemes developed in one-dimensional. There is no this scheme in two-dimensional rigorously established in literature to address sediment transport problems. Here a two-dimensional scheme is developed to address the drawbacks above. The first goal of this paper is to show that the 2D PCCU method accommodates very well to two-dimensional nonconservative sediment transport equations. We will show here that when the conservative laws exist, the 2D PCCU schemes can reproduce some other well-known schemes such as classical 2D path-conservative, 2D CU schemes, 2D path-conservative HLL schemes, 2D HLL solver, and so on. Out of cell these schemes, the PCCU has seen the least far interest in more general two-dimensional nonlinear hyperbolic of nonconservative systems related to sediment transport. This numerical method has originally developed for shallow water equations by Castro et al., **[22]** and was recently extended for Saint-Venant-Exner with a novel well-balanced discretization strategy by Ngatcha et al., **[23]**. The PCCU scheme has the advantage to combine conservative and nonconservative terms discretization and can achieve a high order of accuracy easily with the use of high-order polynomials reconstruction. The presence of sediment transportation/deposition, sediment exchange, friction terms, and bedload equation modifies the design of the scheme. Some numerical methods lose when the sediment transport and morphodynamics are investigated. For sediment transport problems, the numerical scheme must ensure the C-property, captures the shocks and preserves the positivity of water depth. A two-dimensional strategy of well-balanced discretization is developed here to capture the steady-state solutions. We develop here, a 2D hydrostatic reconstruction that preserves positive water depth for all reconstructed values.

*Numerical strategies and Flux approximation techniques*

The numerical strategy used here to solve the STM proposed here is the coupled numerical method. In this strategy, the fluid model, sediment concentration, and morphodynamic model are solved at the same time. The interest of this method is that all the unknowns of the system are updated at the same time steps during the simulation. The discrete flux can be evaluated using three techniques. The first consists to consider only the flux on the center edges of each cell (see Fig.1b). The interest of this strategy is that we get rid of empirical considerations on the evaluation of the flux. The second technique consists to calculate the flux on the vertex and edges of each cell (see Fig.1c). The third strategy consists to calculate the flux only at the vertex of each cell (see Fig.1a). In this paper, the fluxes are evaluated only on the edges of each cell (see Fig. 1b).

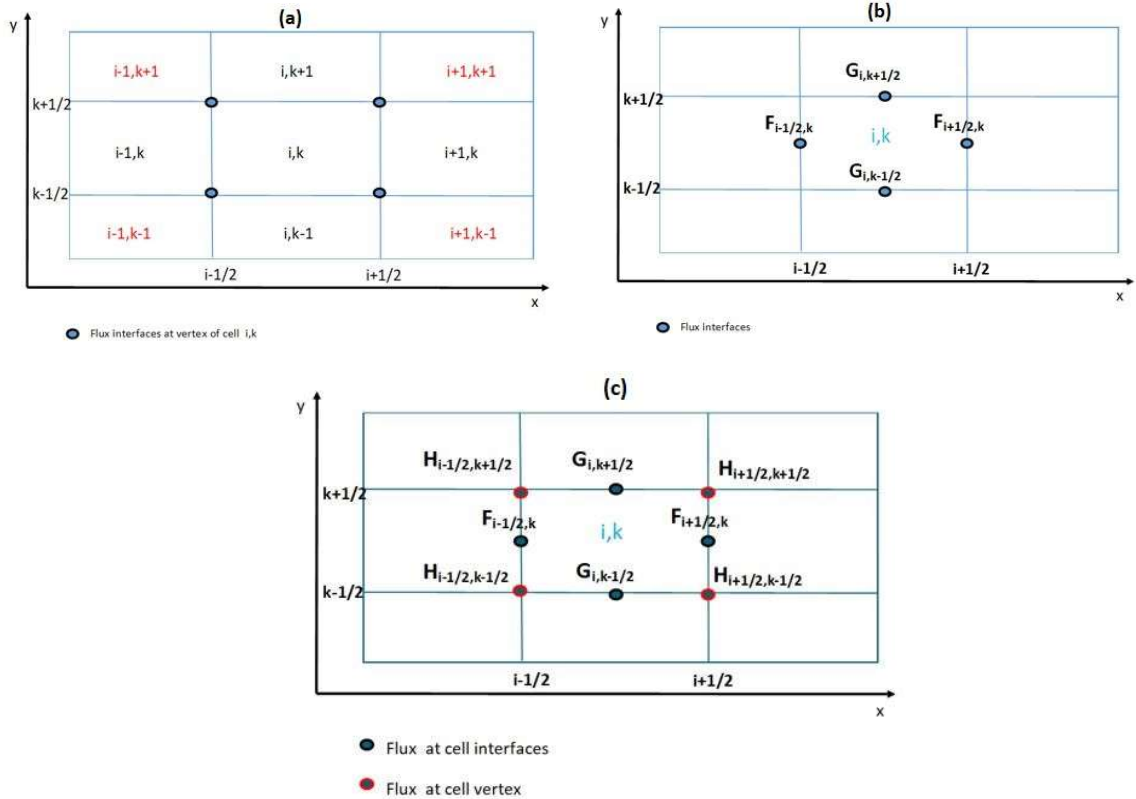

**Figure 1:** 2D Finite volume gridding. Flux evaluated at the vertex of each cell (a); Flux interfaces at the edges of the cell (b); Total flux contribution cells and vertex of each cell (c).

Here, an STM including arbitrarily sloping sediment beds and associated entropy and energy is proposed and solved by a derived 2D well-balanced preserving-positivity PCCU-AENO method on structured meshes. Moreover, an existence theorem of global weak solutions of the model is established and a convergence study is discussed.

*Objectives of paper*

The main objectives of this paper are to: (i) derive a new sediment transport model in a coastal or estuarine environment. The resulting model can be viewed as the generalization of a class of averaged sediment transport models. (ii) derive first and second-order 2D PCCU schemes on structured meshes. (iii) develop a 2D AENO nonlinear reconstruction technique to achieve the second order of accuracy.

*Goals of paper*

One goal of this paper is to introduce a methodology to design a 2D PCCU-AENO scheme on structured meshes for solving nonconservative equations. Another goal is to propose a physical and mathematical analysis (hyperbolicity, existence theorem and convergence) of the model.

*Highlights of paper*

The highlights of the paper are to: (i) Integrate some physicals and hydro-morphodynamic processes to describe the sediment transport in the coastal environment. (ii) Propose an existence theorem of global weak solutions of the model. (iii) Implementation of a methodology to design 2D second-order structured PCCU scheme.

*Scientific Contributions*

The novelties of this paper are the development of a new bedload model; the existence of global weak solutions and convergence results of the model; the development of 2D PCCU schemes with 2D AENO reconstruction; mathematical and physics analysis of the model; some validations with experimental data.

The rest of the paper is presented as follows. Section (II) is dedicated to introducing the mathematical model which couples generalized Shallow Water equations and sediment transport equations, is presented. We study the hyperbolicity of the model in 1D and 2D cases, we give the Rankine-Hugoniot relations and we study the steady-state solutions of the system. We propose the existence theorem of global weak solutions and we expose a convergence result. In section (III), after a brief preliminary on the path-conservative method, a methodology to design a 2D well-balanced PCCU scheme on structured meshes is developed and exposed. We develop for the first time 2D AENO nonlinear reconstruction to obtain a second-order accuracy. In section (IV) some tests are made and the numerical results are compared and discussed.

## II. MATHEMATICAL MODELLING ET HYPERBOLICITY STUDY
### 1. Governing equations.

First, we consider the two-phase equations where each phase $k = s, f$ (sediment '$s$' or fluid '$f$') satisfies the Navier-Stokes (NS) equations as follows:

$$\frac{\partial \alpha_k \rho_k}{\partial t} + div(\alpha_k \rho_k U_k) = 0,$$

$$\frac{\partial \alpha_k \rho_k \mathbf{u}_k}{\partial t} + \nabla.(\alpha_k \rho_k \mathbf{u}_k \otimes \mathbf{u}_k) + \frac{\partial (\alpha_k \rho_k \mathbf{u}_k w_k)}{\partial z} + \nabla P_k = \mathcal{F}_{k,x,y},$$

$$\frac{\partial \alpha_k \rho_k w_k}{\partial t} + \nabla.(\alpha_k \rho_k \mathbf{u}_k w_k) + \frac{\partial (\alpha_k \rho_k w_k w_k)}{\partial z} + \frac{\partial P_k}{\partial z} = \mathcal{F}_{k,z},$$

$$div(U_f) = 0.$$

(1)

Here, $U_k = (\mathbf{u}_k, w_k), \rho_k, P_k, \alpha_k, \mathcal{F}_k$ are respectively the 3D velocity, the density, the pressure the volume fraction and the source terms of the phase $k = s, f$.

Next, we consider the 3D classical NS equations for the evolution of mixing quantities and sediment volume rate obtained by summing the two phases of the system(1). One has:

$$\frac{\partial \rho}{\partial t} + \frac{\partial (\rho u_i)}{\partial x_i} = 0,$$

$$\frac{\partial \rho u_i}{\partial t} + \frac{\partial (\rho u_i u_j)}{\partial x_j} + \frac{\partial P}{\partial x_i} = \mathcal{F}_i, \qquad (2)$$

$$\frac{\partial u_i}{\partial x_i} = 0,$$

where $u_i$, $i = 1, 2, 3$ is the 3D velocity components,

Hydrostatic assumption gives an analytical formulation of the pressure according to the atmospheric pressure (considered here constant) and the vertical water column.

The hydrostatic assumption consists to neglect the fluid vertical acceleration in the flow i.e. the particular derivation $\frac{d(\rho w)}{dt} = 0$.

This leads to: $P = \int_{Z_b}^{\eta} \rho g dz$,

where $g$ is the gravitational constant, $\eta$ is the free surface and $Z_b$ is the bed interface.

In the pressure term, the mixing density $\rho$ is given by:

$$\rho = \rho_w (1-c) + \rho_s c, \qquad (3)$$

where $c$ is the instantaneous sediment concentration, $\rho_w$, $\rho_s$ are respectively the water density and sediment density (assumed constant in time and space). We consider three layers having different densities: a layer of suspension zone, a layer of clear fluid and a layer of bed-load. The suspension is not potential and can be approximately described by the first equation of the system (2) and by using Fick's law as in [8](see also [24]). These equations describe the evolution of fluid mixing in a domain bounded by a dynamic water surface and water bed. We write the mass balance equations in the Saint-Venant formalism [25]. The conservation of momentum is expressed by using the well-known Newton's second law. This law for a control volume states that the net rate of momentum entering the volume (momentum flux) plus the sum of all external forces acting on the volume be equal to the rate of accumulation of momentum. On the free surface, we consider a no-stress condition. On the sediment bed, we consider a no penetration condition. Considering that, we take into account the kinematic boundary conditions on the moving surfaces. A point on the free surface is

$M(x,y,z,t) = -z + \eta(x,y,t)$, where $\eta$ is a smooth function. Assuming that any particle that is on the free surface at the initial instant will remain so at all instants, we have $\dfrac{dM}{dt} = 0$, where the operator $\dfrac{d(.)}{dt}$ is defined by $\dfrac{d(.)}{dt} = \dfrac{\partial(.)}{\partial t} + (\mathbf{u}.\nabla)(.)$. If the free surface volume exchange per time unit, we have $F_u(t) = -z + \eta(x,y,t)$. At the bottom surface or bed interface, one has $z = Z_b(x,y,t)$. Therefore we can define the kinematic boundary conditions on both moving surfaces.

On the bed surface we have:

$$\frac{\partial}{\partial t} + (\mathbf{u}.\nabla)\eta - u_3(\eta) = \frac{dF_u}{dt}, \tag{4}$$

On the free surface we have:

$$\frac{\partial Z_b}{\partial t} + \mathbf{u}(Z_b)\nabla Z_b = \frac{dF_b}{dt} + u_3(Z_b), \tag{5}$$

where the term $\dfrac{dF_b}{dt}$ with $F_b(t) = x_3(t) - Z_b(t,\mathbf{x}(t))$ describes the erosion/deposition exchange; $\dfrac{dF_u}{dt}$ accounts for the effect of lateral contributions. In this work we have assumed that $\dfrac{dF_u}{dt} = 0$. We neglect also the vertical transport at the bed interface i.e. $u_3(Z_b) = 0$. Here, we consider as respectively the volume of sediment deposited and the volume of sediment eroded on the bed. One has:

$$dF_b = dV_s^D - dV_s^E + \phi^* dF_b, \Rightarrow \frac{dF_b}{dt} = \frac{dV_s^D}{dt} - \frac{dV_s^E}{dt} + \phi^* \frac{dF_b}{dt}, \Rightarrow (1-\phi^*)\frac{dF_b}{dt} = D - E,$$

with $D = \dfrac{dV_s^D}{dt}, E = \dfrac{dV_s^E}{dt}$. To retrieve the generalized shallow water-based equations prop water-based we apply an average along the depth of equations (2) using Leibniz's formula to obtain simplified equations. We take an Eulerian approach for the sediment transport equations, rather than the more computationally expensive Lagrangian approach, and make a macroscopic assumption. We introduce into the model an alternative to the bed load equation. Note that for simplicity the diffusion effects of sediment do not integrate into the bedload equation because the effect of advection is more important than the effect of diffusion near the bed due to turbulence or the presence of strong interaction fluid-fluid. Therefore, the bedload sediment transport must depend on the flow regime and size of the grain and the characteristic velocity of the body sedimentary. These parameters aren't incorporated into the classical Exner equation using a sediment transport empirical formula. We recall that sediment transport formulae predicate sediment transport from a given set of hydrodynamic and physical parameters related to sediment and fluid.

# The model

The final two-dimensional model (also named the alternative formulation the of sediment transport model) developed in this paper, is given by the following system:

$$\frac{\partial h}{\partial t} + \frac{\partial hu}{\partial x} + \frac{\partial hv}{\partial y} = \frac{E-D}{(1-p)}$$

$$\frac{\partial hu}{\partial t} + \frac{\partial}{\partial x}\left(huu + \frac{1}{2}gh^2\right) + \frac{\partial}{\partial y}(huv) + gh\frac{\partial Z_b}{\partial x} + \frac{(\rho_s - \rho_w)}{2\rho}gh\frac{\partial hC}{\partial x} - \frac{(\rho_s - \rho_w)}{2\rho}ghC\frac{\partial h}{\partial x} = -C_f u\|\mathbf{u}\| - \frac{(E-D)}{(1-p)}u(Z_b)$$

$$\frac{\partial hv}{\partial t} + \frac{\partial}{\partial y}(huv) + \frac{\partial}{\partial y}\left(hvv + \frac{1}{2}gh^2\right) + gh\frac{\partial Z_b}{\partial y} + \frac{(\rho_s - \rho_w)}{2\rho}gh\frac{\partial hC}{\partial y} - \frac{(\rho_s - \rho_w)}{2\rho}ghC\frac{\partial h}{\partial y} = -C_f v\|\mathbf{u}\| - \frac{(E-D)}{(1-p)}v(Z_b)$$

$$\frac{\partial(hC)}{\partial t} + \frac{\partial(huC)}{\partial x} + \frac{\partial(hvC)}{\partial y} = \frac{\partial}{\partial x}\left(f_s h \nu_m \frac{\partial C}{\partial x}\right) + \frac{\partial}{\partial y}\left(f_s h \nu_m \frac{\partial C}{\partial y}\right) + (E-D)$$

$$\frac{\partial Z_b}{\partial t} + u(Z_b)\frac{\partial(Z_b)}{\partial x} + v(Z_b)\frac{\partial(Z_b)}{\partial y} = -\frac{E-D}{(1-p)}$$

(6)

Here, $h[m]$ is the water depth, $u, v$ are the averaged $x$-velocity and $y$-velocity respectively (with $\mathbf{u} = (u,v)[m/s]$), $hu$, $hv$ (with $\mathbf{q} = (q_1, q_2) = (hu, hv)[m^2/s]$) are the water discharges in both directions and $y$, $Z_b[m]$ is the bed level. $g$ is the gravitational constant. The friction source term is given by Manning's laws: $C_f = n^2 g h^{-1/3}$, where $n[s/m^{1/3}]$ is the manning coefficient and where $g[m/s^2]$ is the constant gravitational. $p$ is the bed porosity. Here, $\rho_w, \rho_s [Kg/m^3]$, $C[m^3/m^3]$ being the water density, sediment density and sediment concentration volumetric respectively. In momentum equations (suspension zone), for sake simplicity we have taken $\mathbf{u}(Z_b) = \mathbf{u}$.

The transport mode parameter $f_s$ is given by:

$$f_s = \min(1; 2.5e^{-Z}),\tag{7}$$

where $Z = \dfrac{W_s}{\kappa u_*}$ is the Rouse number and where $\kappa$ is von Karman number ($\kappa = 0.4$).

$\mathbf{u}_* = \sqrt{C_f \|\mathbf{u}\|^2}$ is the shear stress velocity.

$E[Kg/m^2/s]$ and $D$ are the erosion and deposition given by [1]

$$E = \begin{cases} \varphi(\theta - \theta_{cr,50})h^{-1}\|\mathbf{u}\|d_{50}^{-0.2}, & \text{if } \theta \geq \theta_{cr,50}, \\ 0, & \text{otherwise;} \end{cases} \quad (8)$$

$$D = W_s(1-C_a)^m C_a$$

The deposition rate of sediments $D$ is almost equal to the vertical flux of particle at the boundary.

For erosion rate $\theta = \dfrac{\mathbf{u}_*}{g(s-1)d_{50}}$ is the Shields parameter and the critical Shields parameter is given by:

$$\theta_{cr,50} = \dfrac{0.3}{(1+1.2D_*)} + 0.055(1 - \exp(-0.02D_*))$$

where $D_*$ is the dimensionless grain size parameter, depending of submerged specific gravity of sediment. $\varphi\left[m^{1.2}\right]$ is a coefficient that controls the erosion force. For sediment deposition, $m$ represents the effect of hindered settling due to high sediment concentration and $W_s$ is the fall velocity of sediment given by:

$$W_s = \sqrt{\left(\dfrac{13.95\nu}{d_{50}}\right)^2 + 1.09 s g d_{50}} - 13.95\dfrac{\nu}{d_{50}},$$

where $\nu$ is the kinematic viscosity of water ($\nu = 1.2 \times 10^{-6}$), $d_{50}$ is the average diameter of sediment particles, $s = \dfrac{\rho_s}{\rho_w} - 1$ is the submerged specific gravity of sediment, where $\rho_s$ is the sediment density and $\rho_w$ the water density. $C_a$ is the local near-bed sediment concentration in volume which can be determined following [26]:

$$C_a = \alpha_c C, \quad (9)$$

where $\alpha_c = \min(2, 1 - p/C)$.

## 2. Non heuristic expression of $\mathbf{u}_b(Z_b)$

Given a sedimentary form $Z_b$, two mechanisms influence its evolution. The advection and the 2D flow variation which act strongly on its evolution on the one hand allowing it to move in the flow and on the other hand modifying its geometry. The term on the left of the bed evolution equation is related to advection which highlights that the movement of sedimentary bodies is influenced by a mechanism directly related to the shape of the bottom. From a general point of view, the dynamics of sediment are a lot shorter than the dynamics of the fluid. The scales of waves ($T_w = 10s$) and tidal ($T_t = 12\ hours$) can be easy to identify. The scale of mean current ($T_c \approx 400s$) is a lot shorter than that tidal. The creation of a dune of sand after a flooding can be an approximation as ($T_{sz} = 3\ days$) the scale of their migration on any distance of $400m$ with a velocity estimated to $2m/day$ can be $T_{sx} = 200\ days$ Then, we have $T_w \ll T_c \ll T_t \ll T_{sz} \ll T_{sx}$. Therefore it is important to differentiate the difference between the sediment velocity and fluid velocity i.e. $\mathbf{u}_b(Z_b) \neq \mathbf{u}_s \neq \mathbf{u}$ (where $\mathbf{u}_s$ is the sediment velocity). The equation is not heuristic. Particularly, it's derived from Shallow Water Exner (SW-Exner) model by assuming the $T_c \ll T_{sz}$. This allows us to consider the hydrodynamic equations as stationary with respect to the bed evolution equation (Exner equation) to find the quasi-stationary solution of the mean current. We consider the following SW-Exner system:

$$\nabla.(h\mathbf{u}) = 0, \text{ in } \Omega$$
$$\nabla.\left(h\mathbf{u}\otimes\mathbf{u} + \frac{1}{2}gh^2\right) + gh\nabla Z_b = -C_f \mathbf{u}\|\mathbf{u}\|, \text{ in } \Omega \quad (10)$$
$$\frac{\partial Z_b}{\partial t} + \frac{1}{(1-p)}\nabla.Q_b = 0 \text{ in } \Omega.$$

Here, the sediment flux at the bed is given by the more general formula following:

$$Q_b = a\mathbf{u}^b,\ (a,b) \in \mathbb{R}^2$$

That integrates a large range of sediment transport flux formulas.

After writing the bed sediment evolution equation in terms of hydrodynamic variables, we integrate the energy equation of stationary model (10):

$$\mathbf{u}\nabla.\left(h\mathbf{u}\otimes\mathbf{u} + \frac{1}{2}gh^2\right) + gh\mathbf{u}\nabla Z_b = -C_f \mathbf{u}\mathbf{u}\|\mathbf{u}\|, \text{ in } \Omega \quad (11)$$

We divide by $gh^3\mathbf{u}$ the equation (11) and we rewrite the bed evolution equation. We find the expression of $\mathbf{u}_b(Z_b)$ by analogy to a sand wave model obtained. Here, we have approximated the bed load sediment flux by $\mathbf{u}_b(Z_b)\nabla Z_b$ instead to $\frac{1}{(1-p)}\nabla.Q_b$ as in the model given by (10). Note that $\mathbf{u}_b(Z_b)$ is the characteristic velocity of the advection of sedimentary body $Z_b$ expressed without calibrated parameter as in Ngatcha et al [8] given by:

$$\mathbf{u}_b(Z_b) = \mathbf{u}_b = \frac{1}{1-p}\frac{\partial Q_b}{\partial h}\left(\frac{1}{1-Fr^2}\right)\frac{\mathbf{u}}{\|\mathbf{u}\|}, \quad (12)$$

where $p$ is the sediment bed porosity, $Fr = \frac{|\mathbf{u}|}{\sqrt{gh}}$ is the Froude number. The parameter

$\dfrac{1}{1-p}\dfrac{\partial Q_b}{\partial h}$ describes the sensibility of the sediment transport to water depth. Therefore, we have still $\dfrac{\partial Q_b}{\partial h} \neq \mathbf{u}$. The characteristic velocity of advection of the quantity $\nabla Z_b$ given by (12) precisely depends on the sensibility of water depth and Froude number. Therefore, the body's sedimentary movement is directed by the flow regime. The proposed model is one of the more general existing in the literature and has some advantages as the capability to integrate several sediment transport flux for $(a,b) \in \mathbb{R}^2$ and to differentiate the water velocity from sediment bed velocity (phase lag). The classical Exner model uses some empiric formulas which give approximations results only. These sediment fluxes formulae are designed on uniform flow assumptions and assume that the sediment velocity is equal to fluid velocity. The above system is given by (6) is proven more appropriate to describe the morphodynamics bed evolutions with accuracy and no requires any empirical consideration.

## 3. Some properties of the model

### 3.1. Rankine-Hugoniot relations

In the following, we will assume that $W_L$, $W_R$ are the left and right states in a Riemann problem. Let us define the average and jump operators by

$[\![.]\!] = (.)_R - (.)_L$ and $\{\{.\}\} = \dfrac{(.)_R + (.)_L}{2}$. The Rankine-Hugoniot relation is given by:

$$[\![hu]\!] = \sigma [\![h]\!],$$
$$\left[\!\!\left[ hu^2 + \dfrac{1}{2}gh^2 \right]\!\!\right] + g\{\{h\}\}[\![Z_b]\!] + \dfrac{g\delta\rho}{2\overline{\rho}}\{\{h\}\}[\![hC]\!] - \dfrac{g\delta\rho}{2\overline{\rho}}\{\{hC\}\}[\![h]\!] = \sigma [\![hu]\!],$$
$$[\![huv]\!] = \sigma [\![hv]\!], \qquad (13)$$
$$[\![F_{corr} huC]\!] = \sigma [\![hC]\!],$$
$$\{\{u_b\}\}[\![Z_b]\!] = \sigma [\![Z_b]\!],$$

where $\sigma$ is jump of discontinuities and where $\overline{\rho} = \{\{\rho\}\}$.

### 3.2. Steady state solutions

One goal of this paper is to study the steady-state solution of this new model which is not trivial. The well-balanced scheme proposed here preserves the 1D steady states ''at lake rest''. Indeed, for a smooth solution, we have the following equations:

$$h \equiv \text{constant}, \ hu \equiv \text{constant in time}, \ Z_b \equiv \text{constant in time}, \ C \equiv \text{constant in time},$$
$$\rho \equiv \text{constant in time},$$
(14)

with the machine accuracy.

The structure of 2D steady-state is not easy, but it is possible to find a quasi 1D steady-state solutions:

$$h \equiv \text{constant}, \ hu \equiv \text{constant}, \ hv \equiv 0, \ \partial_x Z_b \equiv \text{constant in time},$$
$$\partial_y Z_b^* \equiv 0, \ \partial_x C \equiv \text{constant in time}, \ \partial_y C \equiv 0, \ \rho \equiv \text{constant in time},$$
(15)

or

$$h \equiv \text{constant}, \ hv \equiv \text{constant}, \ hu \equiv 0, \ \partial_y Z_b \equiv \text{constant in time},$$
$$\partial_x Z_b \equiv 0, \ \partial_y C \equiv \text{constant in time}, \ \partial_x C \equiv 0, \ \rho \equiv \text{constant in time}.$$
(16)

On the other hand, at a point of discontinuity, the steady solutions should verify the Rankine-Hugoniot jump conditions given by (13) with $\sigma = 0$.

As well as the dissipation entropy is given by

$$\left[\!\left[ \left( g(h+Z_b) + \frac{u^2}{2} + \frac{1}{2}\frac{(\rho-\rho_w)}{\rho} gh \right) hu \right]\!\right] \leq 0$$

The R-H relations allow us to conclude that in the whole domain (where the solution is regular and across the discontinuity), the steady states are preserved. The well-balanced 2D PCCU scheme proposed here respect both the "lake at rest" and "dry-lake". Note that the dry lake is obtained when (14) and (15) reduce to

$$u=0, \ hu=0, \ v=0, \ hv=0$$
(17)

## 4. Hyperbolicity study.

Let us $\mathbf{W} = (h, hu, hv, hC, Z_b)^T$, with $W \in \mathbb{R}^5$ the state vector of conservative variables and $F = (F_1, F_2)$ the physical fluxes. We can rewrite the proposed model Eq. (6) in nonconservative form as follows:

$$\frac{\partial \mathbf{W}}{\partial t} + \frac{\partial F_1(\mathbf{W})}{\partial x} + \frac{\partial F_2(\mathbf{W})}{\partial y} + B_{1x}^*(\mathbf{W})\frac{\partial Z_b}{\partial x} + B_{2x}^*(\mathbf{W})\frac{\partial hC}{\partial x} + B_{3x}^*(\mathbf{W})\frac{\partial h}{\partial x} + \ldots$$
$$\ldots B_{1y}^*(\mathbf{W})\frac{\partial Z_b}{\partial y} + B_{2y}^*(\mathbf{W})\frac{\partial hC}{\partial y} + B_{3y}^*(\mathbf{W})\frac{\partial h}{\partial y} = \hat{\mathbf{S}}(\mathbf{W}); \qquad (18)$$

where $\mathbf{x} = (x, y) \in \Omega \subset \mathbb{R}^2$, $t \in (0, T)$.

The vector unknowns $\mathbf{W}: \mathbb{R}^2 \times \mathbb{R}^+ \to \Upsilon$ is a function from space $(x, y) \in \mathbb{R} \times \mathbb{R}$ and time $t$ to the system's state $\Upsilon$ and each components of the flux $F_1$, $F_2: \Upsilon \to \mathbb{R}^5$ is given by

$$F_1(\mathbf{W}) = \begin{pmatrix} hu \\ huu + \frac{1}{2}gh^2 \\ huv \\ huC \\ 0 \end{pmatrix}, \quad F_2(\mathbf{W}) = \begin{pmatrix} hu \\ huv \\ hvv + \frac{1}{2}gh^2 \\ hvC \\ 0 \end{pmatrix}. \qquad (19)$$

The vectors $B_{1x}^*, B_{2x}^*, B_{3x}^*, B_{1y}^*, B_{2y}^*, B_{3y}^*$ are reads

$$B_{1x}^* = \begin{pmatrix} 0 \\ gh \\ 0 \\ 0 \\ u_b \end{pmatrix}, \quad B_{2x}^* = \begin{pmatrix} 0 \\ \frac{gh\delta\rho}{2\rho} \\ 0 \\ 0 \\ 0 \end{pmatrix}, \quad B_{3x}^* = \begin{pmatrix} 0 \\ -\frac{ghC\delta\rho}{2\rho} \\ 0 \\ 0 \\ 0 \end{pmatrix}, \quad B_{1y}^* = \begin{pmatrix} 0 \\ 0 \\ gh \\ 0 \\ v_b \end{pmatrix}, \quad B_{2y}^* = \begin{pmatrix} 0 \\ 0 \\ \frac{gh\delta\rho}{2\rho} \\ 0 \\ 0 \end{pmatrix}, \quad B_{3y}^* = \begin{pmatrix} 0 \\ 0 \\ -\frac{ghC\delta\rho}{2\rho} \\ 0 \\ 0 \end{pmatrix}$$

(20)

The source term reads:

$$\hat{\mathbf{S}}(\mathbf{W}) = S_e + S_F + S_D,$$

where $S_F, S_e, S_D$ are respectively friction source term, the sediment exchange source term and diffusion source term and given respectively by:

$$S_F = \begin{pmatrix} 0 \\ -C_f u \|\mathbf{u}\| \\ -C_f v \|\mathbf{u}\| \\ 0 \\ 0 \end{pmatrix}, S_e = \begin{pmatrix} \dfrac{E-D}{1-p} \\ -\dfrac{(E-D)u}{(1-p)} \\ -\dfrac{(E-D)v}{(1-p)} \\ E-D \\ -\dfrac{E-D}{1-p} \end{pmatrix}, S_D = \begin{pmatrix} 0 \\ 0 \\ 0 \\ \dfrac{\partial}{\partial x}\left(f_s h \nu_m \dfrac{\partial C}{\partial x}\right) + \dfrac{\partial}{\partial y}\left(f_s h \nu_m \dfrac{\partial C}{\partial y}\right) \\ 0 \end{pmatrix}. \quad (21)$$

The numerical solution of nonconservative problem is completed with boundaries conditions and initial conditions of the form $W = W_0$, on $\mathbb{R}^5 \times \{t = 0\}$.

The system can be written in the form of Eq. and is hyperbolic if the Jacobian matrix defined by Eq. has only real eigenvalues and if a full set of linearly independent eigenvectors exists. Therefore, the Jacobian matrix

$$\mathcal{A}_1(\mathbf{W}) = \begin{pmatrix} 0 & 1 & 0 & 0 & 0 \\ -u^2 + gh - \dfrac{\delta\rho}{2\rho}ghC & 2u & 0 & \dfrac{\delta\rho}{2\rho}gh & gh \\ -uv & v & u & 0 & 0 \\ -uC & C & 0 & u & 0 \\ 0 & 0 & 0 & 0 & u_b \end{pmatrix} \quad (22)$$

The quasi-1D system has five distinct eigenvalues:

$$\lambda_1 = u_b, \ \lambda_{2,3} = u, \ \lambda_{4,5} = u \pm \sqrt{gh} \quad (23)$$

A 2D system is hyperbolic in sense that for each state $\mathbf{W} \in \Omega$ and an outer unitary normal vector $v = (v_1, v_2)$, the matrix given by:

$$\mathcal{A}_v(\mathbf{W}) = \mathcal{A}(\mathbf{W}, v) = v_1 \mathcal{A}_1(\mathbf{W}) + v_2 \mathcal{A}_2(\mathbf{W}) \quad (24)$$

has $N+1$ distinct eigenvalues.

According to equation (24) the two-dimensional system has the following eigenvalues:

$$\lambda_1 = \mathbf{u}_b.v, \ \lambda_{2,3} = \mathbf{u}.v, \ \lambda_{4,5} = \mathbf{u}.v \pm \sqrt{gh}, \quad (25)$$

The eigenvectors for associated eigenvalues for 1D case are given by:

$$E_1 = \begin{pmatrix} 1 \\ -u_b \\ 0 \\ C \\ gh - u^2 + 2uu_b - u_b^2 \end{pmatrix}, E_2 = \begin{pmatrix} \frac{\delta\rho}{2\rho} \\ \frac{\delta\rho}{2\rho}u \\ 0 \\ \frac{\delta\rho}{2\rho}C - 1 \\ 0 \end{pmatrix}, E_3 = \begin{pmatrix} \frac{\delta\rho}{2\rho} \\ \frac{\delta\rho}{2\rho}u \\ 1 \\ \frac{\delta\rho}{2\rho}C - 1 \\ 0 \end{pmatrix}, E_4 = \begin{pmatrix} 1 \\ u - \sqrt{gh} \\ 0 \\ C \\ 0 \end{pmatrix}, E_5 = \begin{pmatrix} 1 \\ u + \sqrt{gh} \\ 0 \\ C \\ 0 \end{pmatrix} \quad (26)$$

The third and fourth eigenvalues correspond to genuinely non-linear characteristic fields in the sense of Lax. While remaining eigenvalues correspond to linearly degenerate characteristic fields. $\lambda_{4,5}$ are associated with shock and rarefaction. Riemann invariants are constant across linearly degenerate waves and rarefaction waves whereas for shock waves generalized jump conditions should be satisfied.

**Remark: A resonance condition**

From the eigenstructure of the proposed model, we can see that the conditions for resonance are satisfied if the free internal wavelength that satisfies the unforced equations coincides with the wavelength of topography forcing. This situation appears in our case when [8]:

$$(u - u_b)^2 = gh, \quad \text{in} \quad \Omega^0 \quad (27)$$

It's convenient to set

$$\mathcal{C} = \left\{ \mathbf{W} \in \Omega^0, \ (u_b - u)^2 = gh \right\}, \quad (28)$$

which is the hypersurface on which all the characteristic fields are linearly degenerated. Therefore, the proposed model can predict bed evolution even in the presence of resonance phenomena. In fact, during evolution, a wavelength can be observed in the bedforms between some distances. We also can observe for some waves, the situation where the flow is near the resonance. In the case of floods with sediment transport, for example, resonance situations could occur only when the flood decelerates slowly. In presence of resonance, the above system can be weakly hyperbolic and in this case the vectors $E_1, E_2, E_3, E_4, E_5$ are linearly dependent.

## 5. An existence theorem of global weak solutions of the model

### 5.1. Definition: Weak solution

Let $\Omega \subset \mathbb{R}^2$ with $\mathrm{x} = (x, y) \in \mathbb{R}^2$ an open domain and let $T > 0$ we consider the system given by (6) with the following initial conditions(IC) :

$$h(\mathbf{x},0) = h_0, \ (hu)(\mathbf{x},0) = h_0 u_0, \ (hv)(\mathbf{x},0) = h_0 v_0, \ (hC)(\mathbf{x},0) = h_0 C_0, \ Z_b(\mathbf{x},0) = Z_{b0}. \quad (29)$$

These IC satisfy the following regularity:

$$h_0 \in L^2(\Omega), \ \sqrt{h_0} \in L^2(\Omega), \ \nabla h_0 \in \left(L^2(\Omega)\right)^2, \ Z_{b0} \in L^2(\Omega), \ \frac{(h_0 u_0)^2}{h_0} \in L^1(\Omega), \ u_{b0} \in L^1(\Omega),$$
$$\nabla Z_{b0} \in \left(L^2(\Omega)\right)^2, \ \nabla C_0 \in \left(L^2(\Omega)\right)^2, \ \nabla \sqrt{h_0} \in \left(L^2(\Omega)\right)^2. \quad (30)$$

We say that $(h, hu, hv, hC, Z_b)$ weak solution of model (6) in $\mathcal{D}'(0,T) = C_c^{\infty}(\Omega \times ]0,T])$ which initial data given by () verifying the entropy inequality for all the test functions $\varphi(x,t)$ ($\varphi$ with compact support) such that $\varphi \in C_c^{\infty}(\Omega \times ]0,T[)$ and $\varphi(x,0) = \varphi_0$, $\varphi(x,T) = 0$, if the model (6) holds in $\mathcal{D}'(0,T)$ and the initial condition holds in $\mathcal{D}'(\Omega)$.

For any $\varphi \in C_c^{\infty}(\Omega \times ]0,T])$ we have for the model the variational formulation following:

$$-h_0 \varphi_0 - \int_{\Omega} h \frac{\partial \varphi}{\partial t} + \int_{\Omega \times ]0,T]} (h\mathbf{u}) \nabla \varphi = \int_{\Omega \times ]0,T]} \frac{E-D}{(1-p)} \varphi,$$

$$-h_0 \mathbf{u}_0 \varphi_0 - \int_{\Omega} h\mathbf{u} \frac{\partial \varphi}{\partial t} + \int_{\Omega \times ]0,T]} (h\mathbf{u} \otimes \mathbf{u}):\nabla \varphi + \int_{\Omega \times ]0,T]} (gh\nabla(Z_b + h))\varphi + \int_{\Omega \times ]0,T]} \left(\frac{(\rho_s - \rho_w)}{2\rho} gh^2 \nabla C\right) \varphi =$$
$$- \int_{\Omega \times ]0,T]} C_f \mathbf{u} \|\mathbf{u}\| \ \varphi - \int_{\Omega \times ]0,T]} \left(\frac{(E-D)}{(1-p)} \mathbf{u}(Z_b)\right) \varphi,$$

$$-h_0 C_0 \varphi_0 - \int_{\Omega} (hC) \frac{\partial \varphi}{\partial t} + \int_{\Omega \times ]0,T]} (h\mathbf{u}C) \nabla \varphi = \int_{\Omega \times ]0,T]} (f_s h v_m) \nabla C \nabla \varphi + \int_{\Omega \times ]0,T]} (E-D) \varphi,$$

$$-Z_{b0} \varphi_0 - \int_{\Omega} Z_b \frac{\partial \varphi}{\partial t} + \int_{\Omega \times ]0,T]} (\mathbf{u}(Z_b) \nabla Z_b) \varphi = -\int_{\Omega \times ]0,T]} \frac{E-D}{(1-p)} \varphi,$$

### 5.2. Energy and entropy relations

We recall that the 2D momentum equations of the proposed model written as follows:

$$\frac{\partial h\mathbf{u}}{\partial t} + \nabla \cdot \left( h\mathbf{u} \otimes \mathbf{u} + \frac{1}{2} gh^2 \right) + gh\nabla Z_b + \frac{(\rho_s - \rho_w)}{2\rho} gh\nabla(hC) - \frac{(\rho_s - \rho_w)}{2\rho} ghC\nabla h = -C_f \mathbf{u} \|\mathbf{u}\| - \frac{(E-D)}{(1-p)} \mathbf{u}(Z_b)$$

An energy equation can be obtained by multiply the 2D momentum equations above by $\mathbf{u}$ and by adding by the mass equation multiply by $gh$. We get the following energy equation

$$\frac{d}{dt}\left[\frac{1}{2}h|\mathbf{u}|^2 + \frac{1}{2}gh^2\right] + gh\mathbf{u}\nabla(Z_b + h) + \frac{(\rho_s - \rho_w)}{2\rho} gh^2 \mathbf{u}\nabla C = -C_f \mathbf{u}^2 \|\mathbf{u}\| - \frac{(E-D)}{(1-p)}(\mathbf{u}^2(Z_b) - gh) \quad (31)$$

In absence of sediment/erosion exchange and friction term, we get:

$$\frac{dE}{dt} + gh\mathbf{u}\nabla(Z_b + h) + \frac{(\rho_s - \rho_w)}{2\rho} gh^2 \mathbf{u}\nabla C = 0 ,  \tag{32}$$

where $E = \left[\frac{1}{2}h|\mathbf{u}|^2 + \frac{1}{2}gh^2 + ghZ_b\right]$ is the mechanical energy of the system.

$$\frac{dE}{dt} + \nabla G \leq 0 , \tag{33}$$

where $G = \left(2g(Z_b + \frac{h}{2}) + \frac{\mathbf{u}}{2} + \frac{(\rho_s - \rho_w)}{2\rho} ghC\right) h\mathbf{u}$

### 5.3. Existence of global weak solution and convergence result.

***Theorem [existence of global weak solutions]***

There exists a global weak solution $(h, hu, hv, hC, Z_b)$ of model given by the system of equations (6) satisfying the energy equality and the entropy inequality (31) and (33) respectively. Moreover, its satisfies also the following inequality:

$$\int_\Omega \frac{d}{dt}\left[h|\mathbf{u}|^2 + \frac{3}{2}gh^2\right] + \int_\Omega gZ_b\partial_t h - \int_\Omega \left(C_f \mathbf{u}^2 \|\mathbf{u}\| - \frac{(E-D)}{(1-p)}\left(\mathbf{u}^2(Z_b) - 2g(h+Z_b)\right) + \left(\frac{(\rho_s - \rho_w)C}{2\rho}\right)S_1\right) \leq 0 \tag{34}$$

***Proposition [error estimates]***

According to relations (31) and (33), the following estimates holds:

$$\|h\|_{L^\infty(0,T;H^1(\Omega))} \leq c; \quad \|\nabla\sqrt{h}\|_{L^2(0,T;H^{-1}(\Omega))} \leq c; \quad \|\nabla hC\|_{L^2(0,T;L^2(\Omega))} \leq c; \quad \|\nabla Z_b\|_{L^\infty(0,T;L^2(\Omega))} \leq c;$$

$$\|\sqrt{h}u\|_{L^\infty(0,T;L^2(\Omega))} \leq c; \quad \|\nabla h\|_{L^\infty(0,T;L^2(\Omega))} \leq c; \quad \|Z_b\|_{L^\infty(0,T;L^2(\Omega))} \leq c; \tag{35}$$

**Convergence results**

We consider here a sequence of approximate global weak solution $(h_n, (hu)_n, (hv)_n, (hC)_n, Z_{b,n})$ such that we have

$$\frac{\partial h_n}{\partial t} + \nabla.(h_n\mathbf{u}_n) = S_1^n ,$$

$$\frac{\partial h_n\mathbf{u}_n}{\partial t} + \nabla.\left(h_n\mathbf{u}_n \otimes \mathbf{u}_n + \frac{1}{2}gh_n^2\right) + gh_n\nabla Z_{b,n} + \frac{(\rho_s - \rho_w)}{2\rho_n}gh_n\nabla(h_nC_n) - \frac{(\rho_s - \rho_w)}{2\rho_n}gh_nC_n\nabla h_n = S_2^n ,$$

$$\frac{\partial(h_nC_n)}{\partial t} + \nabla.(h_n\mathbf{u}_nC_n) = \nabla.(f_sh_n\nu_m\nabla C_n) + S_3^n ,$$

$$\frac{\partial Z_{b,n}}{\partial t} + (\mathbf{u}(Z_b))_n \nabla Z_{b,n} = S_4^n ,$$

where $S_1^n = \dfrac{(E-D)}{(1-p)}$, $S_2^n = -C_{f,n}\mathbf{u}_n\|\mathbf{u}_n\| - \dfrac{(E-D)}{(1-p)}\mathbf{u}_n(Z_b)$, $S_3^n = (1-p)S_1^n$, $S_4^n = -S_1^n$.

We assume that its initial values satisfy (for $c$ constant):

$C_0^n \to C_0$ strongly in $L^2(\Omega)$; $h_0^n \to h_0$ strongly in $L^2(\Omega)$; $Z_{b0}^n \to Z_{b0}$ strongly in $L^2(\Omega)$;

$h_0^n \mathbf{u}_0^n \to h_0 \mathbf{u}_0$ strongly in $L^1(\Omega)$; $\sqrt{h_0^n}\mathbf{u}_0^n \to \sqrt{h_0}\mathbf{u}_0$ strongly in $L^1(\Omega)$;

$h_0^n \mathbf{u}_0^n \mathbf{u}_0^n = \dfrac{(q_0^n)^2}{h_0^n} \to h_0 \mathbf{u}_0 \mathbf{u}_0 = \dfrac{(q_0)^2}{h_0}$ strongly in $L^1(\Omega)$. The following relations holds

$$\|h_0\|_{L^2(\Omega)} \le c; \quad \|\sqrt{h_0}\|_{L^2(\Omega)} \le c; \quad \|h_0\mathbf{u}_0\mathbf{u}_0\|_{L^1(\Omega)} \le c; \quad \|h_0\mathbf{u}_0 C_0\|_{L^2(\Omega)} \le c;$$
$$\|\nabla h_0\|_{(L^2(\Omega))^2} \le c; \quad \|Z_{b0}\|_{L^2(\Omega)} \le c; \quad \|h_0\mathbf{u}_0\|_{L^1(\Omega)} \le c; \quad \|\sqrt{h_0}\mathbf{u}_0\|_{L^1(\Omega)} \le c. \tag{36}$$

Moreover, these values verify the following inequality:

$$\int_\Omega h_0|\mathbf{u}_0|^2 + \varepsilon|h_0|^2 + |h_0 Z_{b0}| \le c, \quad \varepsilon\text{=positive constant} \tag{37}$$

**Theorem [convergence]**

There exists a global weak sequence solution $\left(h_n,(hu)_n,(hv)_n,(hC)_n,Z_{b,n}\right)$ of the system (6) with initial values (36) satisfying (34) and (33).

### III. PATH-CONSERVATIVE BASED METHOD FOR NONCONSERVATIVE EQUATIONS.

This section is devoted to presenting some concepts related to the path-conservative method widely used to solve nonconservative problems of the form (18).

**1. A simple classical path-conservative scheme without any intermediate wave (preliminaries)**

The path-conservative approach is used in this work and especially for non-conservative systems. The main idea of this approach is to split the fluctuation into two paths corresponding to left-moving and right-moving waves arising in the Riemann fan solution. This fluctuation is defined $\forall\ \mathbf{W}^+, \mathbf{W}^- \in \Omega$ as:

$$\mathbf{D}(\mathbf{W}^+,\mathbf{W}^-,\nu) = \int_0^1 \mathcal{A}(\mathbf{\Psi}(s,\mathbf{W}^+,\mathbf{W}^-,\nu))\dfrac{\partial \mathbf{\Psi}(s,\mathbf{W}^+,\mathbf{W}^-,\nu)}{\partial s}ds = \mathbf{D}^-(\mathbf{W}^+,\mathbf{W}^-,\nu) + \mathbf{D}^+(\mathbf{W}^+,\mathbf{W}^-,\nu) \tag{38}$$

where $\nu = (\nu_1,\nu_2)$ is out normal of the edge and where $\mathbf{D}^-(\mathbf{W}^+,\mathbf{W}^-,\nu)$ and $\mathbf{D}^+(\mathbf{W}^+,\mathbf{W}^-,\nu)$ are two continuous functions satisfying the following equation:

$$\mathbf{D}^-(\mathbf{W},\mathbf{W},\nu) = \mathbf{D}^+(\mathbf{W},\mathbf{W},\nu) = 0, \ \mathbf{W} \in \Omega. \tag{39}$$

In Eq. (38), the term $\int_0^1 \mathcal{A}(\Psi(s,\mathbf{W}^+,\mathbf{W}^-,\nu)) \dfrac{\partial \Psi(s,\mathbf{W}^+,\mathbf{W}^-,\nu)}{\partial s} ds$ includes the conservative and nonconservative fluxes and writes as follows:

$$\begin{aligned}
\left[\mathcal{A}(\mathbf{W},\nu)\nabla\mathbf{W}\right]_\Psi &= \int_0^1 \left[\mathcal{A}\left(\Psi(s,\mathbf{W}^+,\mathbf{W}^-),\nu\right) - \sigma I\right] \dfrac{\partial \Psi(s,\mathbf{W}^+,\mathbf{W}^-)}{\partial s} ds \\
&= \int_0^1 \left[A_k\left(\Psi(s,\mathbf{W}^+,\mathbf{W}^-),\nu\right)\right] \dfrac{\partial \Psi(s,\mathbf{W}^+,\mathbf{W}^-)}{\partial s} ds + \sum_{m=1}^{3}\left[B^*_{mk}\left(\mathbf{W}^+,\mathbf{W}^-,\nu\right)\left(\mathbf{W}^+ - \mathbf{W}^-\right)\right] \\
&\quad k=1,2
\end{aligned} \tag{40}$$

The definition of path-conservative schemes strongly depends on the chosen family of paths. A non-optimal choice of path can fail the convergence of the solution and even the well-balanced property of path-conservative schemes.

We have used here the linear paths of the form:

$$\Psi(s,\mathbf{W}^+,\mathbf{W}^-) = s\mathbf{W}^+ + (1-s)\mathbf{W}^-, \ \forall \ \mathbf{W}^+, \mathbf{W}^- \in \Omega. \tag{41}$$

This path corresponds to Volpert's definition of the product. This definition does not integrate any intermediate value or unknowns that lead to solving a nonlinear problem via a Newton method. If the Riemann problem associated contains several intermediate waves, the Newton method can fail due to strong interactions that can appear between these waves. In this new path-conservative-based scheme, any intermediate wave is not considered. Well-balanced path-conservative solution. An approximation of steady states solution or non-trivial steady states can be given for a smooth non trivial stationary solution $\mathbf{W} \in \mathbb{R}^2$, where

$$\widetilde{\mathcal{W}} = \left\{\mathbf{W}(\mathbf{x}), \ \left[\mathcal{A}(\mathbf{W},\nu)\nabla\mathbf{W}\right]_\Psi = 0, \ \nabla\mathbf{W} \neq 0, \ \forall \ \mathbf{x} \in \mathbb{R}^2\right\} \tag{42}$$

The equation $\left[\mathcal{A}(\mathbf{W},\nu)\nabla\mathbf{W}\right]_\Psi = 0$ implies that 0 is eigenvalue of $\mathcal{A}(\mathbf{W},\nu)$ and an associated eigenvector for every $\mathbf{x}$, $\nabla\mathbf{W} \neq 0$. Therefore, given an interval $J \subset \mathbb{R}^2$ such that $\mathbf{W}'(\mathbf{x}) \neq 0$

$$\mathcal{W} = \left\{\mathbf{x} \in J, \ J \subset \mathbb{R}^2, \ \nabla\mathbf{W} \neq 0\right\}. \tag{43}$$

The set $\mathcal{W}$ parametrize all the integral curves where the eigenvalues of $\mathcal{A}(\mathbf{W},\nu)$ are zero.

A curve of $\mathcal{W}$ parametrize an arc of an integral curve of a characteristic field. We can reformulated $\widetilde{\mathcal{W}}$ as follows:

$$\widetilde{\mathcal{W}} = \left\{\mathbf{W}(\mathbf{x}) \in \gamma, \gamma \in \mathcal{W}\right\} \tag{44}$$

Using this parametrization, the numerical scheme for solving the sediment transport problem is said to be well-balanced if the following properties are satisfied:

- The scheme solves exactly any smooth stationary solution $\mathbf{W} \in \widetilde{\mathcal{W}}$,
- The scheme solves up to order $k$ any solution $\mathbf{W} \in \widetilde{\mathcal{W}}$ Note that these properties are strongly connected to the relationship between the paths and $\mathcal{W}$.

- **Finite volume gridding for a path-conservative scheme**

Elementary computational cells centered on $(x_i, y_k) = (i\Delta x, k\Delta y)$, where $i, k \in \mathbb{Z}^2$, are denoted $V_{ik} = [x_{i-1/2}, x_{i+1/2}] \times [y_{k-1/2}, y_{k+1/2}]$, where are the corresponding cell interfaces denoted by half integers. The numerical approximations $\mathbf{W}^{\Delta t}$ is the piecewise constant function such that $\forall\, i, k \in \mathbb{Z}^2, \mathbf{W}^{\Delta t}(t^n, x_i, y_k) = \mathbf{W}^n_{i,k}$ on each cells centered $V_{ik}$ with $t^n = n\Delta t,\ n \in \mathbb{N}$. The initial data of $\mathbf{w}$ is denoted by $\mathbf{W}(0, x_i, y_k) = \mathbf{W}^0_{i,k} \in L^\infty(\mathbb{R}^2)$. Once such grid has been designed, we can define at certain time level $t$ the average value of $\mathbf{w}$ overs $V_{ik}$ as:

$$\overline{\mathbf{W}}_{i,k} = \frac{1}{|V_{ik}|} \int_{V_{ik}} \mathbf{W}(x, y, t) dx dy, \qquad (45)$$

where $mes(V_{ik}) = \Delta x \Delta y$. The set of all the cell on the domain $\Omega$ is denoted by $\mathcal{K}_c$ the subscript 'c' denotes the ''center cell centered''. The set of all the edges of $\mathcal{K}_c$ is denoted $\mathcal{E}_c = \mathcal{E}_c^{ext} \cup \mathcal{E}_c^{int}$ where $\mathcal{E}_c^{ext}$ and $\mathcal{E}_c^{int}$ are respectively the exterior edges and interior edges respectively.

2. **Methodology to design a 2D path-conservative central-upwind (PCCU) scheme on structured mesh.**

   In this section, we propose a strategy to design a two-dimensional version of the PCCU scheme on structured meshes. This scheme is a new path-conservative-based scheme where the conservative flux is evaluated using a central-upwind technique and where the fluctuations are evaluated following the Fig. (1a).

   To this end, we follow a concept developed in [22]. We start by developing a two-dimensional CU scheme in the version of path-conservative using the definition of the path-conservative solution presented above. The semi-discrete two-dimensional CU scheme in

path-conservative form is given by:

$$\begin{aligned}
\frac{d\overline{\mathbf{W}}_{i,k}}{dt} &= -\frac{1}{\Delta x}\left(\mathcal{F}_{i+1/2,k} - \mathcal{F}_{i-1/2,k}\right) - \frac{1}{\Delta y}\left(\mathcal{G}_{i,k+1/2} - \mathcal{G}_{i,k-1/2}\right) + \hat{\mathbf{S}}\left(\overline{\mathbf{W}}_{i,k}\right) \\
&= -\frac{1}{\Delta x}\left(\mathcal{F}_{i+1/2,k} - F_1\left(\mathbf{W}^{-}_{i+1/2,k}\right) - \mathcal{F}_{i-1/2,k} + F_1\left(\mathbf{W}^{+}_{i-1/2,k}\right) + F_1\left(\mathbf{W}^{-}_{i+1/2,k}\right) - F_1\left(\mathbf{W}^{+}_{i-1/2,k}\right)\right) \\
&\quad - \frac{1}{\Delta y}\left(\mathcal{G}_{i,k+1/2} - F_2\left(\mathbf{W}^{-}_{i,k+1/2}\right) - \mathcal{G}_{i,k-1/2} + F_2\left(\mathbf{W}^{+}_{i,k-1/2}\right)\right) + \hat{\mathbf{S}}\left(\overline{\mathbf{W}}_{i,k}\right) \\
&= -\frac{1}{\Delta x}\left(D^{-}_{i+1/2,k} + D^{+}_{i-1/2,k} + F_1\left(\mathbf{W}^{-}_{i+1/2,k}\right) - F_1\left(\mathbf{W}^{+}_{i-1/2,k}\right)\right) \\
&\quad - \frac{1}{\Delta y}\left(D^{-}_{i,k+1/2} + D^{+}_{i,k-1/2} + F_2\left(\mathbf{W}^{-}_{i,k+1/2}\right) - F_2\left(\mathbf{W}^{+}_{i,k-1/2}\right)\right) + \hat{\mathbf{S}}\left(\overline{\mathbf{W}}_{i,k}\right) \\
&= \frac{1}{\Delta x}\left(D^{-}_{i+1/2,k} + D^{+}_{i-1/2,k} + \int_{0}^{1} A_1\left(\mathbf{P}_{i,k}(\mathbf{x})\right)\frac{d\mathbf{P}_{i,k}(\mathbf{x})}{d\mathbf{x}}d\mathbf{x}\right) \\
&\quad - \frac{1}{\Delta y}\left(D^{-}_{i,k+1/2} + D^{+}_{i,k-1/2} + \int_{0}^{1} A_2\left(\mathbf{P}_{i,k}(\mathbf{x})\right)\frac{d\mathbf{P}_{i,k}(\mathbf{x})}{d\mathbf{x}}d\mathbf{x}\right) + \hat{\mathbf{S}},
\end{aligned} \quad (46)$$

where the fluctuations function defined by:

$$D^{\pm}_{i+1/2,k} = \left(\frac{1}{2} \pm \frac{1}{2}\frac{a^{+}_{i+1/2,k} - a^{-}_{i+1/2,k}}{a^{+}_{i+1/2,k} - a^{-}_{i+1/2,k}}\right)\left(F_1(\mathbf{W}^{-}_{i+1/2,k}) - F_1(\mathbf{W}^{+}_{i+1/2,k})\right) \pm \frac{1}{2}\left(\frac{-2a^{+}_{i+1/2,k}a^{-}_{i+1/2,k}}{a^{+}_{i+1/2,k} - a^{-}_{i+1/2,k}}\left(W^{+}_{i+1/2,k} - W^{-}_{i+1/2,k}\right)\right)$$

$$= \frac{1 \pm \lambda_1^{i+1/2,k}}{2}\int_0^1\left[A_1\left(\mathbf{\Psi}(s,\mathbf{W}^+,\mathbf{W}^-),v\right)\right]\frac{\partial\mathbf{\Psi}(s,\mathbf{W}^+,\mathbf{W}^-)}{\partial s}ds \pm \frac{\lambda_0^{i+1/2,k}}{2}\left(W^{+}_{i+1/2,k} - W^{-}_{i+1/2,k}\right)$$

(47)

and with

$$D^{\pm}_{i,k+1/2} = \frac{1 \pm \lambda_1^{i,k+1/2}}{2}\int_0^1\left[A_1\left(\mathbf{\Psi}(s,\mathbf{W}^+,\mathbf{W}^-),v\right)\right]\frac{\partial\mathbf{\Psi}(s,\mathbf{W}^+,\mathbf{W}^-)}{\partial s}ds \pm \frac{\lambda_0^{i,k+1/2}}{2}\left(W^{+}_{i,k+1/2} - W^{-}_{i,k+1/2}\right),$$

(48)

where

$$\lambda_1^{i+1/2,k} = \left(\frac{a^{+}_{i+1/2,k} - a^{-}_{i+1/2,k}}{a^{+}_{i+1/2,k} - a^{-}_{i+1/2,k}}\right), \quad \lambda_0^{i+1/2,k} = \frac{-2a^{+}_{i+1/2,k}a^{-}_{i+1/2,k}}{a^{+}_{i+1/2,k} - a^{-}_{i+1/2,k}},$$

$$\lambda_1^{i,k+1/2} = \left(\frac{b^{+}_{i,k+1/2} - b^{-}_{i,k+1/2}}{b^{+}_{i,k+1/2} - b^{-}_{i,k+1/2}}\right), \quad \lambda_0^{i,k+1/2} = \frac{-2b^{+}_{i,k+1/2}b^{-}_{i,k+1/2}}{b^{+}_{i,k+1/2} - b^{-}_{i,k+1/2}}.$$

(49)

Here, the numerical fluxes $\mathcal{F}_{i+1/2,k}$, $\mathcal{G}_{i,k+1/2}$ are given using a CU technique:

$$\mathcal{F}_{i+1/2,k} = \frac{a^{+}_{i+1/2,k}}{a^{+}_{i+1/2,k} - a^{-}_{i+1/2,k}} F_1(\mathbf{W}^{-}_{i+1/2,k}) - \frac{a^{-}_{i+1/2,k}}{a^{+}_{i+1/2,k} - a^{-}_{i+1/2,k}} F_1(\mathbf{W}^{+}_{i+1/2,k}) - \frac{1}{2}\left( \frac{-2a^{+}_{i+1/2,k} a^{-}_{i+1/2,k}}{a^{+}_{i+1/2,k} - a^{-}_{i+1/2,k}} \left( W^{+}_{i+1/2,k} - W^{-}_{i+1/2,k} \right) \right)$$

$$\mathcal{G}_{i,k+1/2} = \frac{b^{+}_{i,k+1/2}}{b^{+}_{i,k+1/2} - b^{-}_{i,k+1/2}} F_2(\mathbf{W}^{-}_{i,k+1/2}) - \frac{b^{-}_{i,k+1/2}}{b^{+}_{i,k+1/2} - b^{-}_{i,k+1/2}} F_2(\mathbf{W}^{+}_{i,k+1/2}) - \frac{1}{2}\left( \frac{-2b^{+}_{i,k+1/2} b^{-}_{i,k+1/2}}{b^{+}_{i,k+1/2} - b^{-}_{i,k+1/2}} \left( W^{+}_{i,k+1/2} - W^{-}_{i,k+1/2} \right) \right)$$

(50)

The topography source term is discretized using a well-balanced discretization strategy proposed in [4]. In Eq. (46)-(48), the functions $\mathbf{W}^{-} = \mathbf{W}(x_\sigma^{-}) = \mathbf{P}_{i,k}(x_\sigma^{-})$ and $\mathbf{W}^{+} = \mathbf{W}(x_\sigma^{+}) = \mathbf{P}_{i+1,k}(x_\sigma^{+})$ such that $\mathbf{W}^{+} \neq \mathbf{W}^{-}$ with $\lim_{x \to x_\sigma^{-}} \mathbf{W} = \mathbf{W}^{-}$ and $\lim_{x \to x_\sigma^{+}} \mathbf{W} = \mathbf{W}^{+}$, where $x_\sigma$ is a discontinuity point. We denoted $\mathbf{W}^{+}$ and $\mathbf{W}^{-}$ the left and right intermediate values of polynomial reconstruction:

$$\widetilde{\mathbf{W}}(x,y,t) = \sum_{i}\sum_{k} \mathbf{P}_{i,k} \mathcal{X}_{V_{ik}}(x), \quad \mathbf{P}_{i,k} = \left( P^{(1)}_{i,k}, P^{(2)}_{i,k}, ...., P^{(N)}_{i,k} \right)^T \tag{51}$$

Here, $\mathcal{X}$ is the characteristic function, $P_i^{(j)}$ are the polynomials of a certain degree satisfying the conservation and accuracy requirements defined for all $i,k$ by:

$$\frac{1}{|V_{ik}|} \int_{V_{ik}} \mathbf{P}_{i,k}(x,y) dx = \overline{\mathbf{W}}_{i,k}, \text{ and } P^{(j)}_{i,k}(x,y) = W^{(j)}(x,y) + O(|V_{ik}|^s), \quad x,y \in V_{ik}, \tag{52}$$

with $s$ a (formal) order of accuracy. $\mathbf{W}(x,y) = (W^{(1)},....,W^{(N)})^t$ is the exact smooth solution. We are interested in left and right limiting values of reconstruction polynomials, often called boundary extrapolated values. The polynomial reconstruction is used to ameliorate the solution approximations at each mesh $V_{ik}$. The order of the scheme depends on the choice of the $\mathbf{P}_i$ functions. For some smooth solution $\mathbf{W}$, we have:

$$\mathbf{W}^{\pm} = \mathbf{W}(\mathbf{x}_\sigma) + \mathcal{O}(|V_{ik}|^s), \quad \forall (i,k) \in \mathbb{Z}^2 \tag{53}$$

The design of the PCCU scheme requires the choice of sufficiently smooth paths

$$\mathbf{\Psi}_{i+1/2,k}(s) = (\Psi^{(1)}_{i+1/2,k}, \Psi^{(2)}_{i+1/2,k}, ...., \Psi^{(N)}_{i+1/2,k}, \Psi^{(N+1)}_{i+1/2,k}) := \mathbf{\Psi}(s, \mathbf{W}^{-}, \mathbf{W}^{+}) \tag{54}$$

connecting the two states $\mathbf{W}^{-}$ and $\mathbf{W}^{+}$ across the jump discontinuity at $\mathbf{x} = \mathbf{x}_0$ such that a local-Lipschitz application $\mathbf{\Psi}:[0,1]\times\Omega\times\Omega \to \Omega$ satisfies the following property:

$$\mathbf{\Psi}(s, \mathbf{W}^{+}, \mathbf{W}^{-}) = s\mathbf{W}^{+} + (1-s)\mathbf{W}^{-}, \forall \mathbf{W}^{+}, \mathbf{W}^{-} \in \Omega$$
(55)

We have in this scheme taken a simplest linear segment path in each direction:

$$\Psi_{i+1/2,k}(s) = \mathbf{W}^-_{i+1/2,k} + s(\mathbf{W}^+_{i+1/2,k} - \mathbf{W}^-_{i+1/2,k}), \quad \Psi_{i,k+1/2}(s) = \mathbf{W}^-_{i,k+1/2} + s(\mathbf{W}^+_{i,k+1/2} - \mathbf{W}^-_{i,k+1/2}), \quad s \in [0,1]$$
(56)

The values $\mathbf{w}$ at point $(i-1/2,k)$, $(i+1/2,k)$, $(i,k-1/2)$, $(i,k+1/2)$ are given as follows:

$$\begin{aligned}
\mathbf{W}^+_{i+1/2,k} &= \mathbf{P}_{i+1,k}(x_{i+1/2} - 0, y_k), \quad \mathbf{W}^-_{i-1/2,k} = \mathbf{P}_{i,k}(x_{i-1/2} - 0, y_k) \\
\mathbf{W}^+_{k+1/2,i} &= \mathbf{P}_{i,k+1}(x_i, y_{k+1/2} - 0), \quad \mathbf{W}^-_{k-1/2,i} = \mathbf{P}_{i,k}(x_i, y_{k-1/2} + 0)
\end{aligned}$$
(57)

Note that all the above quantities depend on time, but we simplify the notation by suppressing this dependence.

**Proposition III.1**

The one-sided local speeds of propagation $a^\pm_{i+1/2,k}$ and $b^\pm_{i+1/2,k}$ are upper/lower bounds on the largest/smallest eigenvalues of Jacobian matrix given above:

$$\begin{aligned}
a^+_{i+1/2,k} &= \max\left\{u^-_{i+1/2,k} + \sqrt{gh^-_{i+1/2,k}}, u^+_{i+1/2,k} + \sqrt{gh^+_{i+1/2,k}}, u^+_{i+1/2,k}, u^-_{i+1/2,k}, u^+_{b,i+1/2,k}, u^-_{b,i+1/2,k}, 0\right\}; \\
a^-_{i+1/2,k} &= \min\left\{u^-_{i+1/2,k} - \sqrt{gh^-_{i+1/2,k}}, u^+_{i+1/2,k} - \sqrt{gh^+_{i+1/2,k}}, u^+_{i+1/2,k}, u^-_{i+1/2,k}, u^+_{b,i+1/2,k}, u^-_{b,i+1/2,k}, 0\right\}; \\
b^+_{i,k+1/2} &= \max\left\{u^-_{i,k+1/2} + \sqrt{gh^-_{i,k+1/2}}, u^+_{i,k+1/2} + \sqrt{gh^+_{i,k+1/2}}, u^+_{i,k+1/2}, u^-_{i,k+1/2}, u^+_{b,i,k+1/2}, u^-_{b,i,k+1/2}, 0\right\}; \\
b^-_{i,k+1/2} &= \min\left\{u^-_{i,k+1/2} - \sqrt{gh^-_{i,k+1/2}}, u^+_{i,k+1/2} - \sqrt{gh^+_{i,k+1/2}}, u^+_{i,k+1/2}, u^-_{i,k+1/2}, u^+_{b,i,k+1/2}, u^-_{b,i,k+1/2}, 0\right\}.
\end{aligned}$$
(58)

Moreover, the CFL condition reads:

$$\Delta t \leq CFL \min\left(\frac{\Delta x}{4a}, \frac{\Delta y}{4b}\right); \quad 0 < CFL < 1, \quad a = \max(a^+_{i+1/2,k}, -a^-_{i+1/2,k}), \quad b = \max(b^+_{i,k+1/2}, -b^-_{i,k+1/2}),$$
(59)

where $\Delta t$ is the step time.

**Remark III.2**

For conservative equations, the fluctuation terms given by equations (47) and (48) contains only the terms associated to derivative of flux:

$$\int_0^1 A_\zeta(\mathbf{P}_{i,k}(\mathbf{x})) \frac{d\mathbf{P}_{i,k}(\mathbf{x})}{d\mathbf{x}} d\mathbf{x} = F_\zeta(\mathbf{W}^-_{i+1/2,k}) - F_\zeta(\mathbf{W}^+_{i-1/2,k}), \text{ and}$$

$$\int_0^1 \left[A_\zeta(\Psi(s, \mathbf{W}^+, \mathbf{W}^-), v)\right] \frac{\partial \Psi(s, \mathbf{W}^+, \mathbf{W}^-)}{\partial s} ds = F_\zeta(\mathbf{W}^+_{i+1/2,k}) - F_\zeta(\mathbf{W}^-_{i+1/2,k}), \quad \zeta = 1, 2. \quad (60)$$

For nonconservative systems, we will write $\mathcal{A}_\zeta$ given above instead of $A_\zeta$. When the fluxes are computed in CU sense, the resulting path-conservative central-upwind scheme is a version of path-conservative HLL Riemann solver.

## 3. The 2D PCCU scheme on structured meshes

For two-dimensional path-conservative central-upwind method without topography source term, the fluctuation terms are now given by:

$$\mathbf{D}_{i,k+1/2}^{\pm} = \frac{1 \pm \lambda_1^{i,k+1/2}}{2} \int_0^1 \left[ \mathcal{A}_2\left(\boldsymbol{\Psi}(s, \mathbf{W}^+, \mathbf{W}^-), \nu\right) \right] \frac{\partial \boldsymbol{\Psi}(s, \mathbf{W}^+, \mathbf{W}^-)}{\partial s} ds \pm \frac{\lambda_0^{i,k+1/2}}{2} \left( \mathbf{W}_{i,k+1/2}^+ - \mathbf{W}_{i,k+1/2}^- \right) \quad (61)$$

$$\mathbf{D}_{i+1/2,k}^{\pm} = \frac{1 \pm \lambda_1^{i+1/2,k}}{2} \int_0^1 \left[ \mathcal{A}_1\left(\boldsymbol{\Psi}(s, \mathbf{W}^+, \mathbf{W}^-), \nu\right) \right] \frac{\partial \boldsymbol{\Psi}(s, \mathbf{W}^+, \mathbf{W}^-)}{\partial s} ds \pm \frac{\lambda_0^{i+1/2,k}}{2} \left( \mathbf{W}_{i+1/2,k}^+ - \mathbf{W}_{i+1/2,k}^- \right) \quad (62)$$

According to the scheme given by the equation (46), a two-dimensional version of the PCCU scheme can be easily designed on structured meshes.

The definition of fluctuation given by (61) and (62) permit to show that the Central-upwind scheme can be written as in version of path-conservative without major difficulty.

The first order semi-discrete PCCU scheme reads:

$$\frac{d\overline{\mathbf{W}}_{i,k}}{dt} = -\frac{1}{\Delta x}\left(\mathcal{F}_{i+1/2,k} - \mathcal{F}_{i-1/2,k}\right) - \frac{1}{\Delta y}\left(\mathcal{G}_{i,k+1/2} - \mathcal{G}_{i,k-1/2}\right) - \left(\sum_{m=1}^3 \left(B_m^*\right)^{\boldsymbol{\Psi}}\right)_{i,k}$$
$$+ \frac{1}{\Delta x \Delta y}\left[\frac{a_{i+1/2,k}^-}{a_{i+1/2,k}^+ - a_{i+1/2,k}^-}\left(\sum_{m=1}^3 \left(B_{mx}^*\right)^{\boldsymbol{\Psi}}_{i+1/2,k}\right) - \frac{a_{i+1/2,k}^+}{a_{i+1/2,k}^+ - a_{i+1/2,k}^-}\left(\sum_{m=1}^3 \left(B_{mx}^*\right)^{\boldsymbol{\Psi}}_{i-1/2,k}\right)\right] \quad (63)$$
$$\frac{1}{\Delta x \Delta y}\left[\frac{b_{i,k+1/2}^-}{b_{i,k+1/2}^+ - b_{i,k+1/2}^-}\left(\sum_{m=1}^3 \left(B_{my}^*\right)^{\boldsymbol{\Psi}}_{i,k+1/2}\right) - \frac{b_{i,k+1/2}^+}{b_{i,k+1/2}^+ - b_{i,k+1/2}^-}\left(\sum_{m=1}^3 \left(B_{my}^*\right)^{\boldsymbol{\Psi}}_{i,k-1/2}\right)\right] + \hat{S}_{i,k}$$

where we have let

$$\left(B_{mx,y}^*\right)_{i,k} = \int_{V_{i,k}} \left(B_m^*\right)\left(\mathbf{P}_{i,k}(\mathbf{x})\right) \left(\frac{dP_{i,k}^{(1)}}{d\mathbf{x}}, \frac{dP_{i,k}^{(2)}}{d\mathbf{x}}, \ldots, \frac{dP_{i,k}^{(N)}}{d\mathbf{x}}\right)^T d\mathbf{x}, \quad (64)$$

$$\left(B_{mx}^*\right)^{\boldsymbol{\Psi}}_{i+1/2,k} = \int_0^1 \left(B_{mx}^*\right)\left(\boldsymbol{\Psi}_{i+1/2,k}(s)\right) \left(\frac{d\psi_{i+1/2,k}^{(1)}}{ds}, \ldots, \frac{d\psi_{i+1/2,k}^{(N)}}{ds}\right)^T ds, \quad m=1,2,3 \text{ and}$$

$$\left(B_{my}^*\right)^{\boldsymbol{\Psi}}_{i,k+1/2} = \int_0^1 \left(B_{my}^*\right)\left(\boldsymbol{\Psi}_{i,k+1/2}(s)\right) \left(\frac{d\psi_{i,k+1/2}^{(1)}}{ds}, \ldots, \frac{d\psi_{i,k+1/2}^{(N)}}{ds}\right)^T ds, \quad m=1,2,3 \quad (65)$$

Using the linear path, a very accurate numerical approximation of the characteristic velocity of body sedimentary also can be given by:

$$u_b^* = \int_0^s u_b(s)ds = \sum_{g=1}^{NGp} w_g u_b(s_g), \tag{66}$$

where $NGp$ is a number of points Gauss quadrature rule, $w_g$ are the weights and $s_g$ are the positions distributed in the unit interval $[0, 1]$. Here we have considered:

$$s_1 = \frac{1}{2}, \ s_{2,3} = \frac{1}{2} \pm \frac{\sqrt{15}}{10}, \ w_1 = \frac{8}{18}, \ w_{2,3} = \frac{5}{18}. \tag{67}$$

In all the numerical simulation one point-Gauss quadrature is used and therefore we have

$$u_b^* = \frac{8}{18} u_b\left(\frac{1}{2}\right). \tag{68}$$

This choice allows us to ensure the achievement of second of accuracy.

The semi-discrete first order PCCU scheme is given by Equations (50)-(58) and (63)-(65). To achieve the second-order, we use an AENO-type reconstruction technique in conjunction with ADER schemes for hyperbolic equations.

## 4. 2D AENO nonlinear reconstruction procedure and properties of the scheme

### *4.1 2D AENO reconstruction*

Here, we describe a new second-order extension of the PCCU in space using a modified version of the averaging essentially non-oscillatory (AENO) procedure originally developed in one-dimensional by Toro et al., [27]. Here, an original two-dimensional version of AENO nonlinear reconstruction is developed to improve numerical solutions (which allows achieving the second-order accuracy in space of the scheme). We start by writing a 2D piecewise operator of the form:

$$\mathbf{P}_{i,k}(x,y) = \overline{\mathbf{W}}_{i,k} + \Delta^x_i(x - x_i) + \Delta^y_k(y - y_k); \ x, y \in V_{ik}, \ x_i = \frac{x_{i+1/2} - x_{i-1/2}}{2}, \ y_k = \frac{y_{k+1/2} - y_{k-1/2}}{2}$$
(69)

where $\Delta_{i,k} = (\nabla \mathbf{W})_{i,k}$ are the slopes that approximate $(\nabla \mathbf{W}(x_i, y_k, t^n))$ in a non-oscillatory manner using a nonlinear slope obtained by convex combination of $\Delta^x_i$ and $\Delta^y_k$ as follows

$$(\nabla \mathbf{W})^x_{i,k} = \Delta^x_{i,k} = \beta^x \Delta^n_{i+1/2,k} + (1-\beta^x)\Delta^n_{i-1/2,k}, \ \beta^x \in [0,1] \tag{70}$$

$$(\nabla \mathbf{W})^y_{i,k} = \Delta^y_{i,k} = \beta^y \Delta^n_{i,k+1/2} + (1-\beta^y)\Delta^n_{i,k-1/2}, \ \beta^y \in [0,1] \tag{71}$$

where

$$\beta^x(r^x) = \frac{r^x}{\sqrt{l^2 + r^{x2}}}; \text{ with } r^x = \frac{|\Delta_{i-1/2,k}|}{\Delta_{i+1/2,k} + \epsilon}, \quad \beta^y(r^y) = \frac{r^y}{\sqrt{l^2 + (r^y)^2}}; \text{ with } r^y = \frac{|\Delta_{i,k-1/2}|}{\Delta_{i,k+1/2} + \epsilon} \text{ and }$$

where

$$\Delta_{i+1/2,k} = \frac{\overline{W}_{i+1,k} - \overline{W}_{i,k}}{\Delta x}, \quad \Delta_{i-1/2,k} = \frac{\overline{W}_{i,k} - \overline{W}_{i-1,k}}{\Delta x}, \quad \Delta_{k+1/2,i} = \frac{\overline{W}_{i,k+1} - \overline{W}_{i,k}}{\Delta y}, \quad \Delta_{k-1/2,i} = \frac{\overline{W}_{i,k} - \overline{W}_{i,k-1}}{\Delta y},$$

$l$ is a positive parameter, $\epsilon$ is a small positive tolerance to avoid division by zeros.

The resulting semi-discrete second-order two-dimensional PCCU-AENO scheme for the is then given by Equations (50)-(58) and (63)-(69).

### 4.2. Well-balanced property

At the steady states, $\forall i,k$

$$\overline{h}_{i,k}^n = \overline{h}_{i+1/2,k}^{n,\pm} = h_0, \ \overline{hC}_{i,k}^n = \overline{hC}_{i+1/2,k}^{n,\pm} = K_0, \ \overline{Z}_{b,i,k}^n = \overline{Z}_{b,i+1/2,k}^{n,\pm} = b_0, \ \overline{u}_{i,k}^n = \overline{u}_{i+1/2,k}^{n,\pm} = 0, \ E - D = 0, \ \rho_{i,k}^n = C_1$$
(72)

Thus, $\overline{h}_{i,k}^n, \overline{(hu)}_{i,k}^n, \overline{(hv)}_{i,k}^n, \overline{hC}_{i,k}^n, \overline{Z}_{b,i,k}^n$ are reconstructed, for all steps time.

Then, we have also:

$$\overline{\eta}_{i,k}^n = \overline{\eta}_{i+1/2,k}^{n,\pm} = \overline{h}_{i,k}^n + \overline{Z}_{b,i,k}^n = \overline{h}_{i+1/2,k}^{n,\pm} + \overline{Z}_{b,i+1/2,k}^{n,\pm} = \eta_1 \tag{73}$$

i.e. $\overline{\eta}_{i,k}^n$ is constant at the lake at rest steady states. According to the above equations, we have:

$$\mathbf{W}_{i,k+1/2}^+ - \mathbf{W}_{i,k+1/2}^- = 0 \text{ and } \mathbf{W}_{i+1/2,k}^+ - \mathbf{W}_{i+1/2,k}^- = 0 \tag{74}$$

At steady states we have:

$$\begin{aligned}
\mathcal{F}_{i+1/2,k}^{(1)} - \mathcal{F}_{i-1/2,k}^{(1)} + S_{e,i,k}^{(1)} &= 0 \\
\mathcal{F}_{i+1/2,k}^{(2)} - \mathcal{F}_{i-1/2,k}^{(2)} - B_{1x,i,k}^{*(2)} - B_{2x,i,k}^{*(2)} + \mathcal{H}_{i+1/2,k}^{(2)} - \mathcal{H}_{i-1/2,k}^{(2)} &= 0 \\
\mathcal{F}_{i+1/2,k}^{(3)} - \mathcal{F}_{i-1/2,k}^{(3)} + S_{e,i,k}^{(3)} &= 0 \\
\mathcal{F}_{i+1/2,k}^{(4)} - \mathcal{F}_{i-1/2,k}^{(4)} + S_{e,i,k}^{(4)} &= 0
\end{aligned} \tag{75}$$

We have therefore the following well balanced discretization topography term:

$$B_{1x,i,k}^{*(2)} + B_{2x,i,k}^{*(2)} + B_{3x,i,k}^{*(2)} = \mathcal{F}_{i+1/2,k}^{(2)} - \mathcal{F}_{i-1/2,k}^{(2)} + \mathcal{H}_{i+1/2,k}^{(2)} - \mathcal{H}_{i-1/2,k}^{(2)} \tag{76}$$

Here, the numerical flux $\mathcal{F}_{i+1/2}^{(2)}$ is given in CU sense and reads:

$$\mathcal{F}^{(2)}_{i+1/2,k} = \frac{a^{+}_{i+1/2,k}}{a^{+}_{i+1/2,k} - a^{-}_{i+1/2,k}} (0.5g(h^{-}_{i+1/2,k})^2) - \frac{a^{-}_{i+1/2,k}}{a^{+}_{i+1/2,k} - a^{-}_{i+1/2,k}} (0.5g(h^{+}_{i+1/2,k})^2). \tag{77}$$

The nonconservative contribution $\mathcal{H}^{(2)}_{i+1/2,k}$ reads

$$\mathcal{H}^{(2)}_{i+1/2,k} = \frac{a^{-}_{i+1/2,k}}{a^{+}_{i+1/2,k} - a^{-}_{i+1/2,k}} \left( B^{\Psi(2)}_{1x,i+1/2,k} + B^{\Psi(2)}_{2x,i+1/2,k} + B^{\Psi(2)}_{3x,i+1/2,k} \right) \tag{78}$$

where $\quad B^{\Psi(2)}_{i+1/2,k} = -g \dfrac{\left( h^{+}_{i+1/2,k} + h^{-}_{i+1/2,k} \right)}{2} \left( h^{+}_{i+1/2,k} - h^{-}_{i+1/2,k} \right)$  (79)

The well-balanced scheme is obtained by replacing the discrete topography term $B^{\Psi(2)}_{i+1/2,k}$ given in (65) by that given by Eq. (79).

- ***Well-balanced discrete source terms***

Here, we use the reconstructed unknowns to discretize the source term in well-balanced sense. The terms $S_e, S_D$ and $S_F$ are discretized as follows:

$$S_{e,i,k} = \begin{pmatrix} \dfrac{E-D}{1-p} \\ -\dfrac{(E-D)}{(1-p)} u_{i,k} \\ -\dfrac{(E-D)}{(1-p)} v_{i,k} \\ E-D \\ -\dfrac{E-D}{1-p} \end{pmatrix} \text{ and } S_{F,i,k} = \begin{pmatrix} 0 \\ -g \dfrac{\left( h^{-}_{i+1/2,k} + h^{+}_{i-1/2,k} \right)}{2} S_{fx,i,k} \\ -g \dfrac{\left( h^{-}_{i,k+1/2} + h^{+}_{i,k-1/2} \right)}{2} S_{fy,i,k} \\ 0 \\ 0 \end{pmatrix} \tag{80}$$

$$S_{D,i,k} = \begin{pmatrix} 0 \\ 0 \\ 0 \\ \left( f_{s,i+1/2,k} h_{i+1/2,k} v_f \dfrac{C_{i+1,k} - C_{i,k}}{dx} \right) + \left( f_{s,i,k+1/2} h_{i,k+1/2} v_f \dfrac{C_{i,k} - C_{i,k-1}}{dy} \right) \\ 0 \end{pmatrix}$$

With (80), the proposed 2D AENO-PCCU scheme satisfies the C-property.

### 4.3. Preserving-positivity reconstruction

Here we expose a discretization strategy developed in [23] that preserves positive the water depth. This procedure has been improved for a one-dimensional total sediment transport in [8]. The 2D version of this methodology is presented in follows:

We note $\mathbf{W} = (\eta = h + Z_b, hu, hv, hC, Z_b)$, where $\eta$ is the free surface. Let us

$$\mathcal{W}_\tau = \left\{ \mathbf{W}_{i,k}^n = \left( \overline{h}_{i,k}^n, \overline{q}_{1i,k}^n, \overline{q}_{2i,k}^n, \overline{hC}_{i,k}^n, \overline{Z}_{b,i,k}^n \right) \in \mathbb{R}_\tau^5, \ \overline{h}_{i,k}^n \succ \overline{h}_{i\pm1/2,k}^{n,+}, \overline{h}_{i,k\pm1/2}^{n,-} \geq 0 \right\}, \tag{81}$$

the discrete admissible space that preserves the positivity of water depth. The left/right velocities and concentrations are calculated as:

$$u_{i+1/2,k}^\pm = \frac{(hu)_{i+1/2,k}^\pm}{h_{i+1/2,k}^\pm}, \quad v_{i+1/2,k}^\pm = \frac{(hv)_{i+1/2,k}^\pm}{h_{i+1/2,k}^\pm} \quad \text{and} \quad C_{i+1/2,k}^\pm = \frac{(hC)_{i+1/2,k}^\pm}{h_{i+1/2,k}^\pm}, \tag{82}$$

The bottom reconstruction at left and right is given by:

$$\overline{Z}_{b,i+1/2,k}^+ = \min\left( \max\left( \overline{Z}_{b,i,k}, \overline{Z}_{b,i+1,k} \right), \overline{\eta}_{i+1,k} \right), \quad \overline{Z}_{b,i+1/2,k}^- = \min\left( \max\left( \overline{Z}_{b,i,k}, \overline{Z}_{b,i+1,k} \right), \overline{\eta}_{i,k} \right), \tag{83}$$

where $\overline{\eta}_{i+1,k} = \overline{h}_{i+1,k} + \overline{Z}_{b,i+1,k}$, $\overline{\eta}_{i,k} = \overline{h}_{i,k} + \overline{Z}_{b,i,k}$.

Which verified at the steady states:

$$\overline{h}_{i+1/2,k}^{*,-} = \overline{h}_{i+1/2,k}^{*,+}, \quad \overline{h}_{i+1/2,k}^- = \overline{h}_{i+1/2,k}^+, \tag{84}$$

where

$$\overline{h}_{i+1/2,k}^{*,-} = \min\left( \overline{\eta}_{i,k} - \max\left( \overline{Z}_{b,i+1/2,k}^-, \overline{Z}_{b,i+1/2,k}^+ \right), \overline{h}_{i,k} \right), \quad \overline{h}_{i+1/2,k}^{*,+} = \min\left( \overline{\eta}_{i+1,k} - \max\left( \overline{Z}_{b,i+1/2,k}^-, \overline{Z}_{b,i+1/2,k}^+ \right), \overline{h}_{i+1,k} \right)$$

Therefore, the preserving-positivity reconstruction of water depth is given by:

$$h_{i+1/2,k}^\pm = \max\left( 0, \ h_{i+1/2,k}^{*,\pm} \right) \tag{85}$$

Now, the steady states are verified by:

$$\overline{h}_{i+1/2,k}^{*,-} = \overline{h}_{i+1/2,k}^{*,+}, \quad \overline{h}_{i+1/2,k}^- = \overline{h}_{i+1/2,k}^+. \tag{86}$$

Then, we recomputed the averaged variables $\overline{q}_1, \overline{q}_2$ and $\overline{hC}$ at each interface as

$$\overline{q}_{1,i+1/2,k}^\pm = \overline{h}_{i+1/2,k}^\pm \overline{u}_{i+1/2,k}^\pm, \quad \overline{q}_{2,i+1/2,k}^\pm = \overline{h}_{i+1/2,k}^\pm \overline{v}_{i+1/2,k}^\pm \quad \text{and} \quad \left( \overline{hC} \right)_{i+1/2,k}^\pm = \overline{h}_{i+1/2,k}^\pm \overline{C}_{i+1/2,k}^\pm. \tag{87}$$

We use the same methodology in y-direction without any difficulty to obtain the steady state equations following:

$$\overline{h}_{i,k+1/2}^{*,-} = \overline{h}_{i,k+1/2}^{*,+}, \quad \overline{h}_{i,k+1/2}^- = \overline{h}_{i,k+1/2}^+. \tag{88}$$

Therefore, we have $\bar{h}_i^{n+1} \in \mathcal{W}_\tau$ since for all $i$, $\bar{h}_i^n \in \mathcal{W}_\tau$ at time $t = t^n$.

**Remark**

- The proposed 2D PCCU-AENO scheme can be applied to arbitrary nonconservative systems.
- The resulting scheme can used to solve three-dimensional hyperbolic problems.
- It proven that the proposed fully discrete 2D PCCU-AENO scheme is well balanced and preserves the positivity of water depth.

## IV. RESULTS AND DISCUSSION

In this section, several 1D and 2D tests are performed to assess the performances of the proposed 2D PCCU method and model. The Fully discrete scheme computed is simply obtained via a strong stability preserving (SSP) based method used in [6]. For some tests, the error estimate between the numerical solution and the reference solution is computed and the convergence rate is deduced. More generally, the error estimate is evaluated in $L^1 -norm$ at the time $t = T$, where $T$ is the final time. For one-dimensional tests, the numerical stability is imposed by the Courant-Friedrich-Lewy (CFL) condition. The integration time step is evaluated as:

$$\Delta t = CFL \frac{mes(K_i)}{2a}, \qquad 0 < CFL \leq 1$$

where $K_i = [x_{j-1/2}, x_{j-1/2}]$ and where $a = \max(a^-_{j+1/2}, -a^-_{j+1/2})$, $a^\pm_{j+1/2}$ being the local propagation speeds. AENO reconstruction is performed using $\epsilon = 0.0001$, $l = 1$.

The computational parameters for some tests are given in Table 1.

**Table 1:** Parameter values for the simulation

| Parameter | $\rho_w$ | $\rho_s$ | $\varphi$ | $\nu$ | $p$ | $g$ | $d_{50}$ | $n$ | $m$ |
|---|---|---|---|---|---|---|---|---|---|
| Value | 1000 | 2650 | 0.015 | 0.000012 | 0.4 | 9.8 | 0.001 | 0.028 | 2 |

### 1. Verification of C-property.

This test is designed to verify that, when the erosion/deposition exchange source term is zero, all the rest of contributions will not affect the well-balanced property of the scheme. Such test is also done in [4]. We show that our scheme is able to exactly preserve the steady states at rest. The initial conditions are:

$$\begin{aligned} &h(x,y,0) = 2 - Z_b(x,y,0), \text{ with } Z_b(x,y,0) = 0.02 + 0.1\exp((-x-0.5)^2 - (y-0.5)^2) \\ &\text{and } u(x,y,0) = 0, v(x,y,0) = 0, C(x,y,0) = 0.7\exp(-5(x-0.9)^2 - 50(y-0.5)^2) \end{aligned} \qquad (89)$$

We uses zero-order extrapolation at all of the boundaries. The initial condition is displayed in Fig. 2. The domain of simulation is $\Omega = [0,1] \times [0,1]$. We run the 2D PCCU method with AENO reconstruction using 400 structured cells and the obtained results are displayed in Fig. 3. at time $t = 10s$. The results shows that our well balanced discretization of bed slope terms preserves exactly "lake at rest" which is still physically significant.

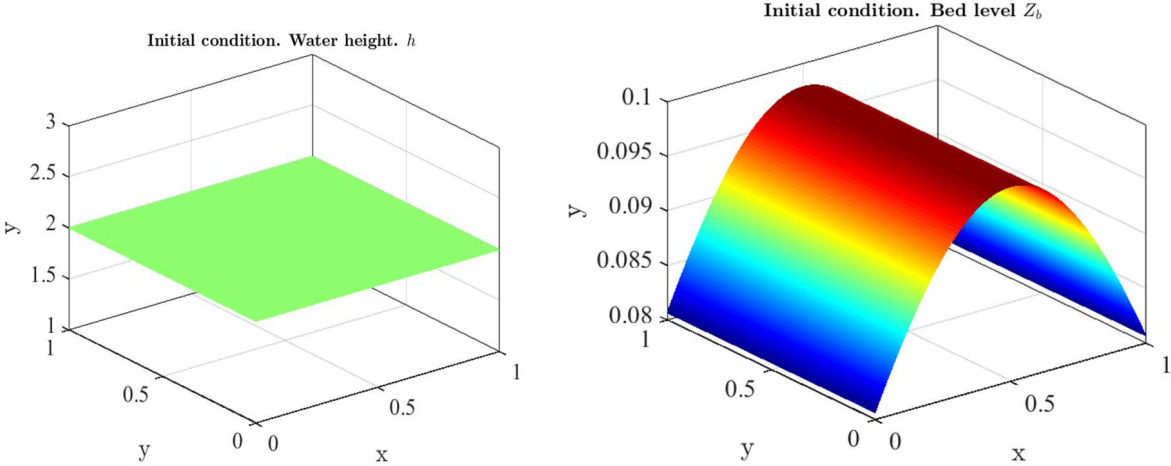

**Figure 2:** Initial condition for well-balanced test.

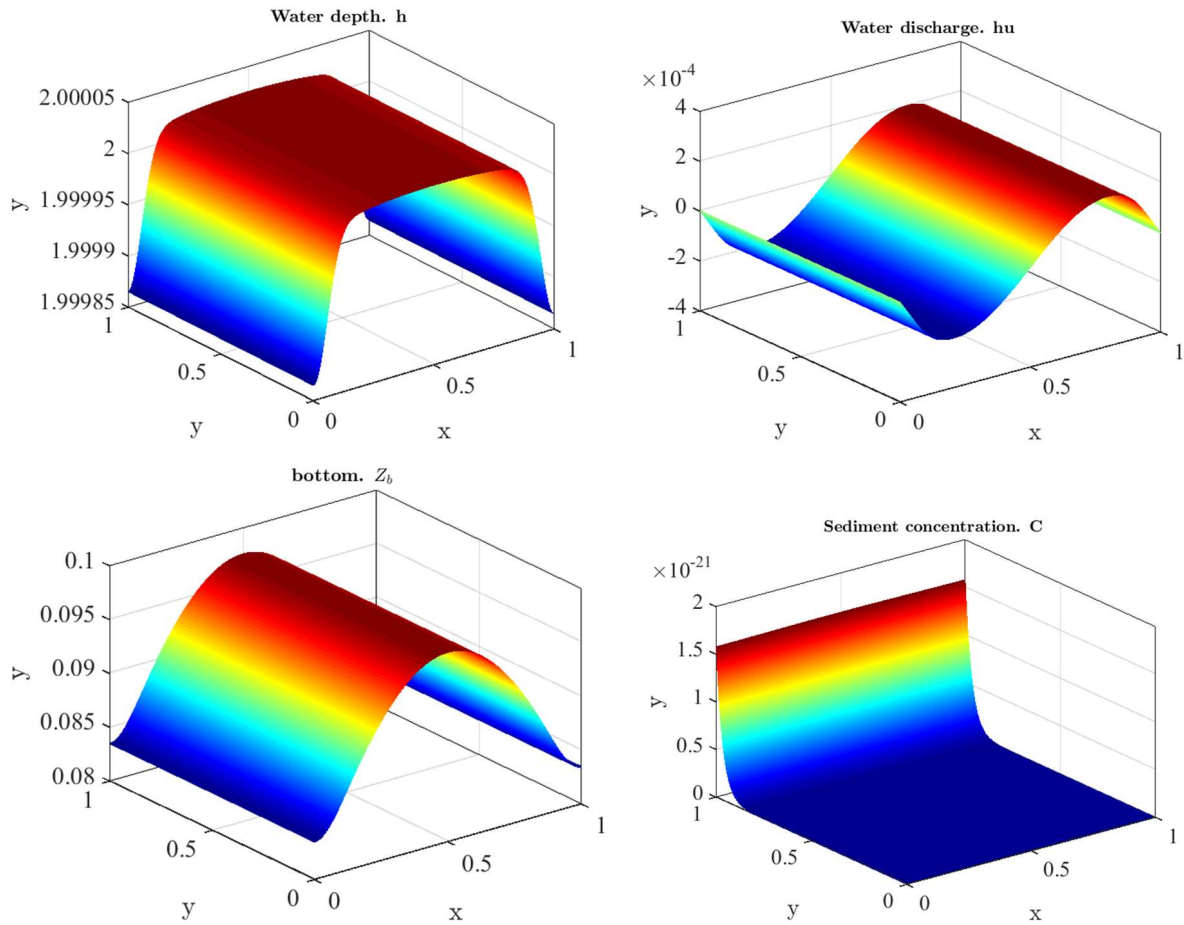

**Figure 3:** Computational solutions of well-balanced test. Water height $h$, bed level $Z_b$, water discharge $hu$ and deposit mass $C$ profiles at time t=10s.

Really, in the nature the situation where the water not moving in a river or channel is impossible. In this test, very small variations of the sediment bed, velocity, sediment concentration and water depth are observed during the simulation. This is correct according to the observation in the nature. The constant water depth cannot be observed in reality. For sediment concentration, a small variation is observed at the beginning of the simulation after that the stable equilibria are retrieved. It is expected that the water-free surface remains practically constant and the sediment concentration should be zero at all times. The sediment concentration is practically zeros during the simulation since the variation scale is very negligible compared to the rest. The water discharge of fluid remains constant during these exchanges. The water discharge varies around zero. All these small variations are observed on a microscopic scale to see really the behavior of steady-state solutions during a long time simulation. Therefore, the proposed PCCU-AENO method preserves the C-property to the machine's precision. We verify the convergence of the proposed method by using the measure of the difference between the solutions computed on two consecutive grids. The $L^1-norm$ is given by:

$$\|\Phi^N - \Psi^N\|_1 = \frac{1}{N^2} \sum_{i=1}^{N} \sum_{k=1}^{N} |\Phi_{i,k}^N - \Psi_{i,k}^N|, \tag{90}$$

where $\Phi^N := \{\Phi_{i,k}^N\}$ and $\Psi^N := \{\Psi_{i,k}^N\}$ are two functions prescribed on structured mesh of $N \times N$ cells. The rates of convergence are calculated as:

$$\mathcal{O}(L^1) = Log_2 \left( \frac{\|\varphi^{N/2} - \varphi^{N/4}\|_1}{\|\varphi^N - \varphi^{N/2}\|_1} \right), \tag{91}$$

where we have noted that $Log_b(x) = y \Rightarrow b^y = x$.

**Table 2:** Estimate error for well-balanced test.

|  | h |  | hu |  | hC |  | $Z_b$ |  |
|---|---|---|---|---|---|---|---|---|
| N | $L^1$ | $\mathcal{O}(L^1)$ | $L^1$ | $\mathcal{O}(L^1)$ | $L^1$ | $\mathcal{O}(L^1)$ | $L^1$ | $\mathcal{O}(L^1)$ |
| 400 | 4.034E-4 | / | 2.87E-2 | / | 7.343E-4 | / | 2.044E-4 | / |
| 800 | 1.018E-3 | 1.98 | 6.348E-3 | 2.07 | 1.547E-4 | 1.96 | 1.31E-3 | 1.92 |
| 1600 | 2.448E-4 | 2.06 | 1.708E-3 | 2.06 | 4.001E-5 | 2.05 | 3.41E-4 | 1.79 |
| 3200 | 6.082E-5 | 1.99 | 4.01E-4 | 1.95 | 9.457E-6 | 2.01 | 9.08E-3 | 2.03 |

## 2. Experimental validation 1D test.

In this test, the 1D version of the model is solved and the results are compared with experimental data and classical Exner model. A similar test is done by [11] using explicit staggered finite volume scheme and by using 1D PCCU scheme. We test the capability of our model to reproduce the sediment transport even in an experimental channel. The initial conditions are given by:

$$h(x,0) = \begin{cases} 0.1 & \text{if } x \leq 0 \\ 0 & \text{if } x > 0 \end{cases}, \quad u(x,0) = 0, \quad Z_b(x,0) = 0, \quad E - D = 0. \tag{92}$$

For the classical Exner model, the sediment diameter is $d_{50} = 0.0032$, sediment density is $\rho_s = 1.540$, the domain of simulation is $\Omega = [-1.25; 1.25]$. Grass formula is used for $Q_b$. The free surface $\eta = h + Z_b$ and bed level profiles at different times $t = 0.5$, $t = 0.7$, $t = 1$ using our proposed model are shown in Fig. 4, and those obtained using classical Shallow Water Exner model are plotted in Fig. 5. They show a good agreement between the numerical computation and the experimental data (available in [28] see also [29]) with respect to the water level and sediment profiles. We have used in all the simulation $CFL = 0.1$, $N = 100$ cells.

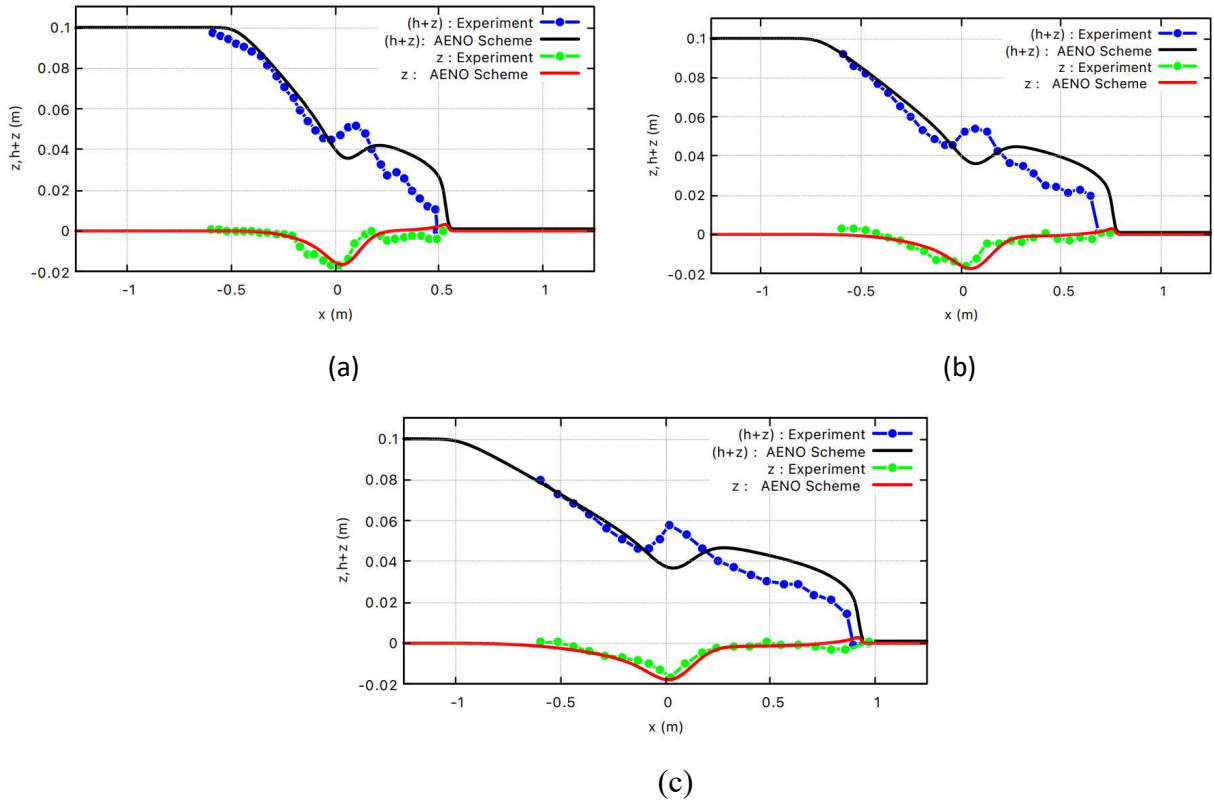

**Figure 4:** Computational solution of the proposed model using PCCU. Comparison with experimental data.

We observe that the water level and sediment bed profiles are better approximated using PCCU-AENO scheme. The waves of the model are well captured during the simulation. These profiles are different from those obtained by using the classical Shallow water Exner model as presented in Fig. (5).

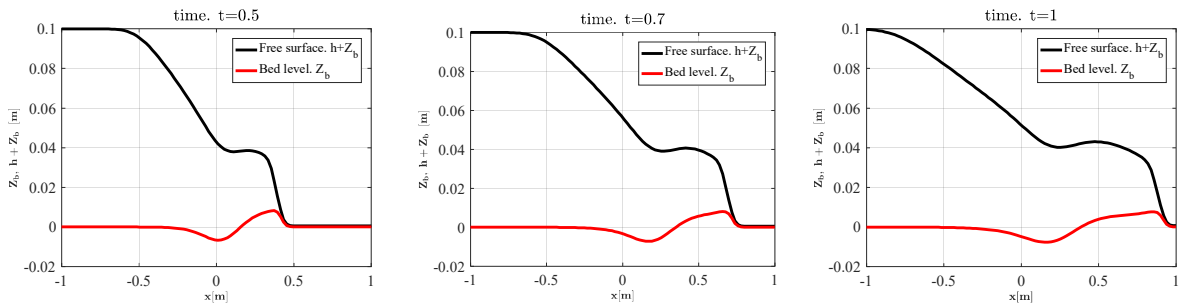

**Figure 5:** Computational solutions obtained by Shallow Water Exner model with Grass formula, $CFL=0.5$, $N=100$.

The classical Exner model coupled with a bed-load sediment flux formula widely used to describe the morphodynamics of coastal environments does not give good results according to experimental data. However, the main drawbacks of this model remains its lack of robustness. This observation is also done in [8].

## 3. Multiple grains size test. Sediment diffusion effect.

We perform now the same previous test with erosion/deposition effect $E - D \neq 0$ and with sediment diffusion effect. We use the same initial conditions as in the previous test (experimental validation test). Ones compare different profiles of sediment concentration using PCCU-AENO scheme with different sediment diameters $d_1 = 0.002, d_1 = 0.0032, d_1 = 0.008, d_1 = 0.02 (mm)$. It's well known that the deposition/erosion exchange depends strongly of sediment diameter (see the formula of these function in appendix). The obtained results are plotted in Fig. (6). The test shows that the proposed model is able to simulate a wide range of sediment class size.

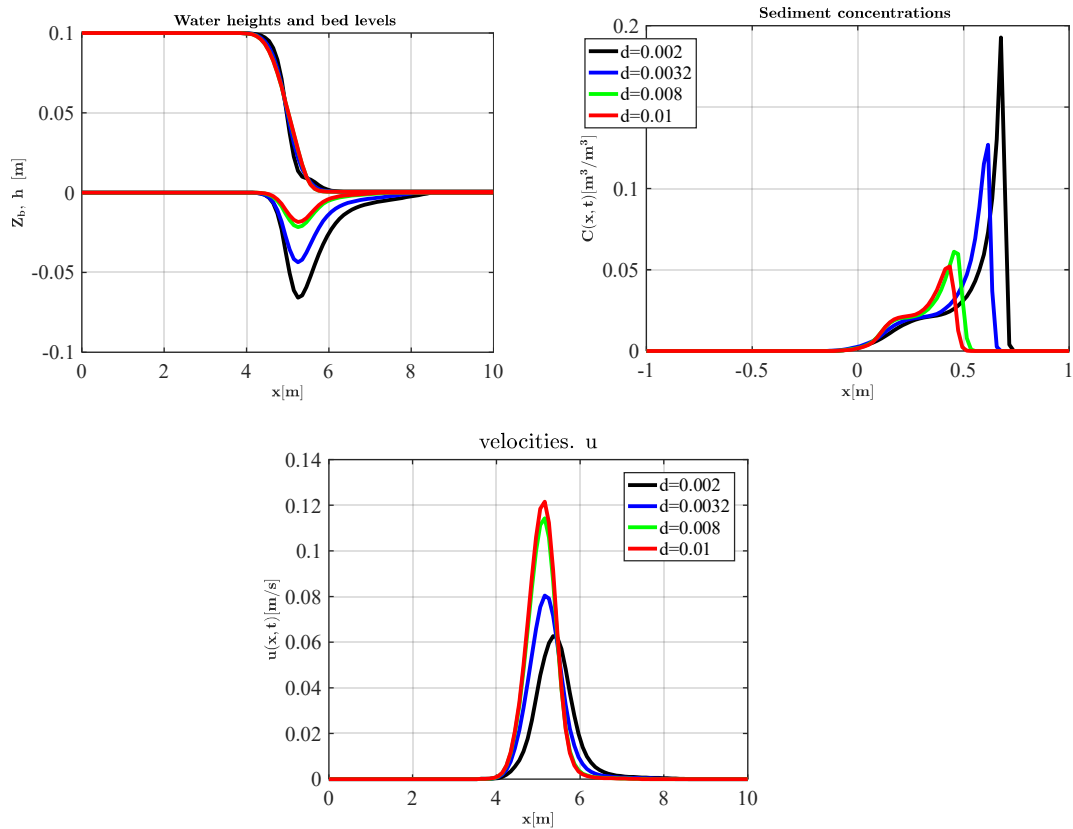

**Figure 6:** Numerical solution of sediment concentration, bed level, water height and velocities using PCCU-AENO scheme. Comparison between different sediment diameters. We have used N=100 cells, t=0.25, $CFL = 0.1$.

We expected that our bed sediment model does not depend on sediment diameter as the classical Exner model. Classical models use some empiric formula that gives approximate results only on a range of flow regimes and sediment diameters. Some of these formulas become uncertain when the sediment diameter becomes greater. The diffusion effect is well visible in the profiles of sediment concentration. It's observed that the sediment concentration is more adapted for fine grains which are associated with low velocity due to mixing flow. The presence of sediment in the water reduces its flow velocity. When the size of the sediment becomes greater, the concentration becomes low and fluid/sediment velocity has the same behavior as fluid velocity. The profiles obtained by our simulation are in agreement with what could be observed in nature or an experimental channel. Particularly, the profiles of sediment concentration are interesting and very close to the results obtained by [30] even if they are not the same conditions.  The proposed shock-capturing scheme can serve to produce more realistic simulations in real environment conditions.

## 4. Bed evolution movement.

We study here the bed evolution when the sediment bed is not fixed. The initial conditions are given by:

$$h(x,y,0) = 1 - Z_b(x,y,0), \quad \text{with} \quad Z_b(x,y,0) = 0.02 + 0.1\exp((-x-0.5)^2 - (y-0.5)^2)$$
$$\text{and} \quad u(x,y,0) = 0, v(x,y,0) = 0, \ C(x,y,0) = 0.01. \tag{93}$$

This initial values are displayed in Fig. (7). The numerical solution obtained by applied 2D well-balanced positivity-preserving PCCU scheme is plotted in Fig. (8).

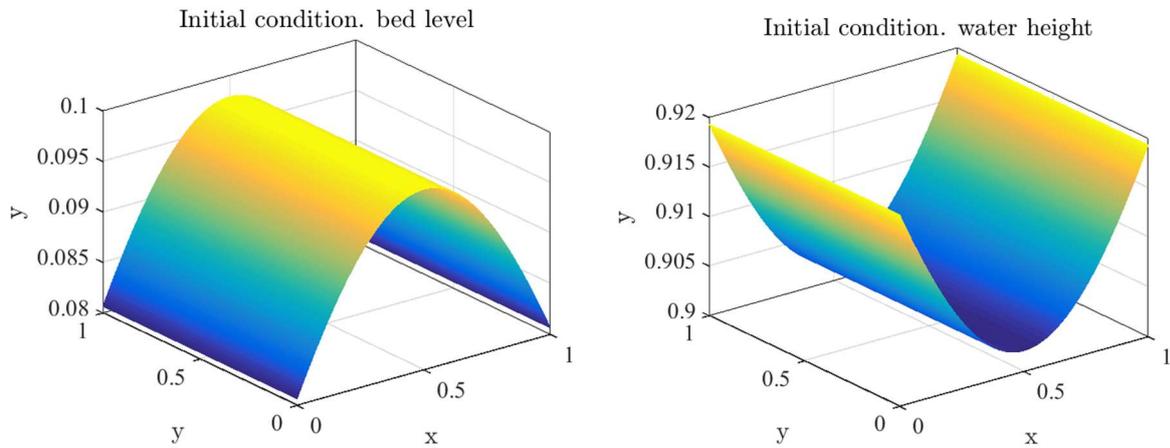

**Figure: 7**   Initial condition. Bed and water height profiles.

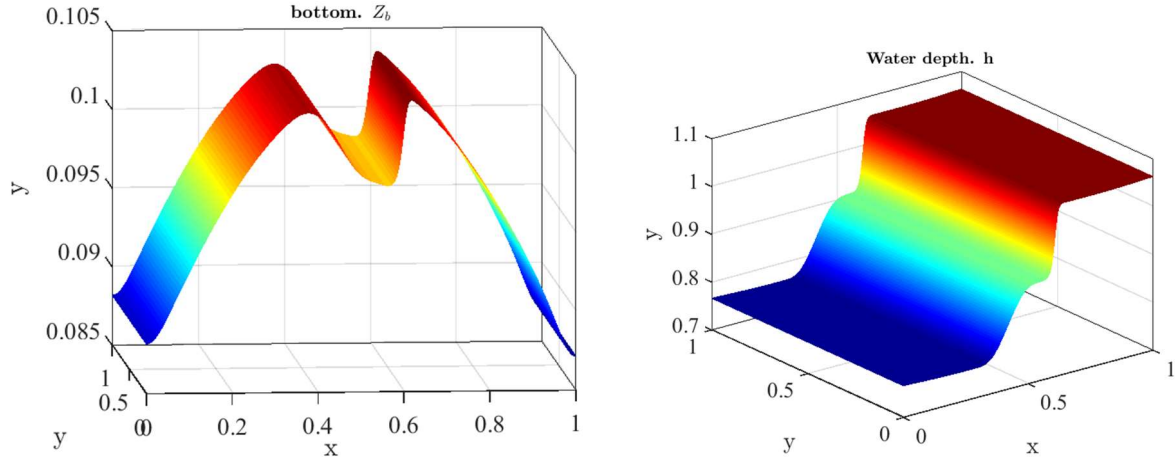

**Figure: 8** Movement of sediment bed and free surface at $t = 0.3s$ with $N_x = N_y = 400$ cells.

The movement of bed is well described and the water height profile is agreement with respect to the physic of the problem studied here. The movement of the bed and the water level are well computed and the phase lag effect is well accounts. These profiles are well observed in the nature.

## 5. 2D Riemann problem

We consider here, the 2D Riemann problem with initial data given in Table. (2). This Riemann problem consists of dam-break over erodible bed with sediment transport. We recall that the initial condition for the local Riemann problem is given by:

$$\mathbf{W}(x,y,0) = \begin{cases} \mathbf{W}_{LW} & \text{if } x<0, y>0 \\ \mathbf{W}_{RW} & \text{if } x<0, y>0 \\ \mathbf{W}_{LD} & \text{if } x<0, y>0 \\ \mathbf{W}_{RD} & \text{if } x<0, y>0 \end{cases}$$

This test simulates rapid spatial and temporal 2D deformations of the free surface and sediment bed. The boundary condition is free, the computational domain is $\Omega = [-1,1] \times [-1,1]$.

**Table. 3** Initial condition for the Riemann problem

| Domain | $h(x,y,0)$ | $u(x,y,0)$ | $v(x,y,0)$ | $Z_b(x,y,0)$ |
|---|---|---|---|---|
| $x \in [-0.5, 0.5]; y \in [-0.5, 0.5]$ | 2 | 0 | 0 | 2 |
| $x \in [-1,1]; y \in [-1,1]$ | 1 | 0 | 0 | 1 |

The initial concentration volume is $C = 0.001$. The rest of computational parameters is given by Table. 1. We can analyze the solution of this two-dimensional Riemann problem computed on uniform grid $N_x = 400$, $N_y = 400$ cells. The initial conditions is given in Fig. (9).

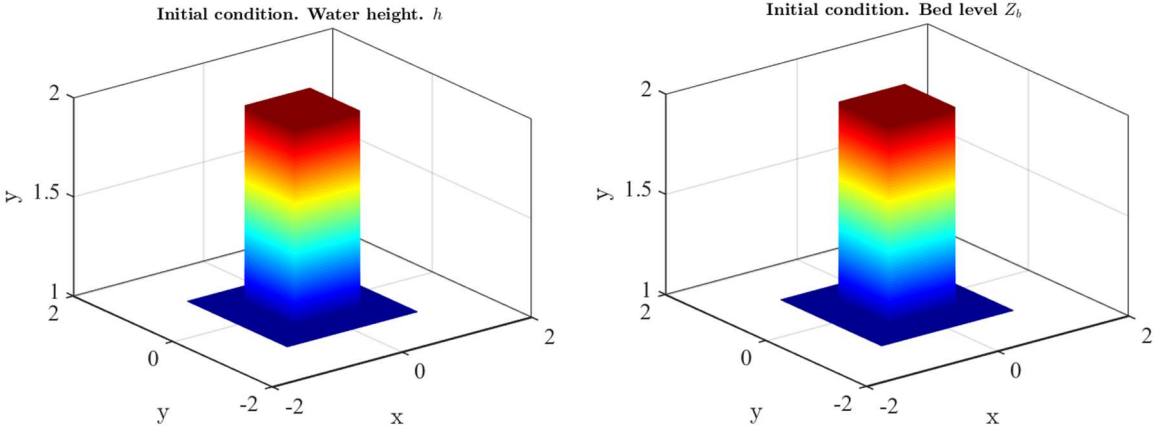

**Figure: 9** Initial conditions for 2D Riemann problem

The computed solution of the Riemann problem using 2D PCCU scheme are plotted in Fig. (10).

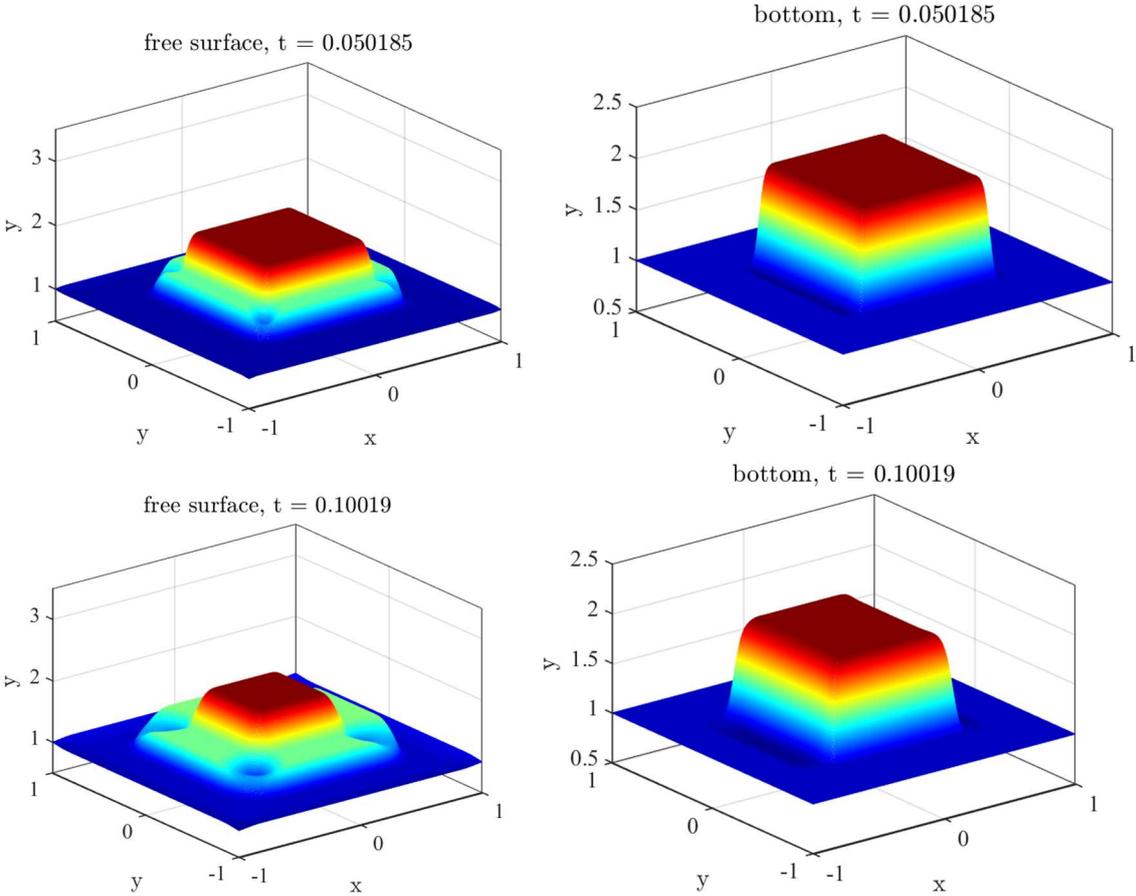

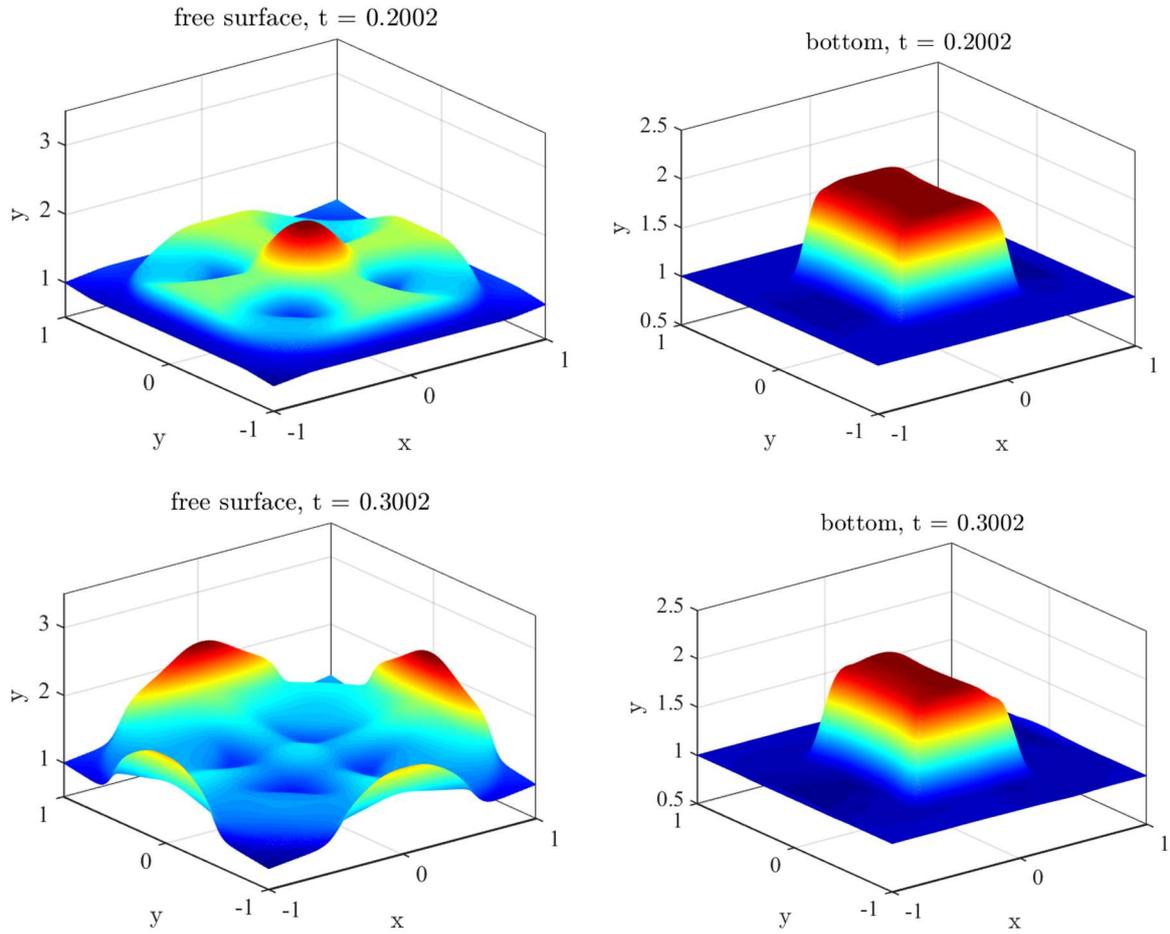

**Figure: 10** Computational solution of the Riemann problem. Bed level and free surface profiles after fourth simulations. CFL=0.5.

We plot here the sediment concentration and the bed evolution profiles when the time is in Fig. (11). It's observed interesting physic related to the dynamic of sediment. We expected that after a long time, the sediment deposition/erosion exchange is very high. The 2D profile of sediment concentration is well captured during the simulation. It's very rare to find a such test in the literature showing really a 2D behavior of sediment transport in regular channels.

**Table 4**  Estimate error for 2D Riemann problem.

|   | $h$ |   | $Z_b$ |   |
|---|---|---|---|---|
| $N$ | $L^1$ | $\mathcal{O}(L^1)$ | $L^1$ | $\mathcal{O}(L^1)$ |
| 100 | 9.818E-3 | / | 2.11E-2 | / |
| 200 | 2.717E-3 | 1.73 | 5.83E-3 | 1.85 |
| 400 | 7.188E-4 | 1.961 | 1.34E-3 | 1.972 |
| 800 | 1.887E-4 | 1.922 | 3.58E-4 | 1.882 |
| 1600 | 4.571E-5 | 1.990 | 8.91E-5 | 1.992 |

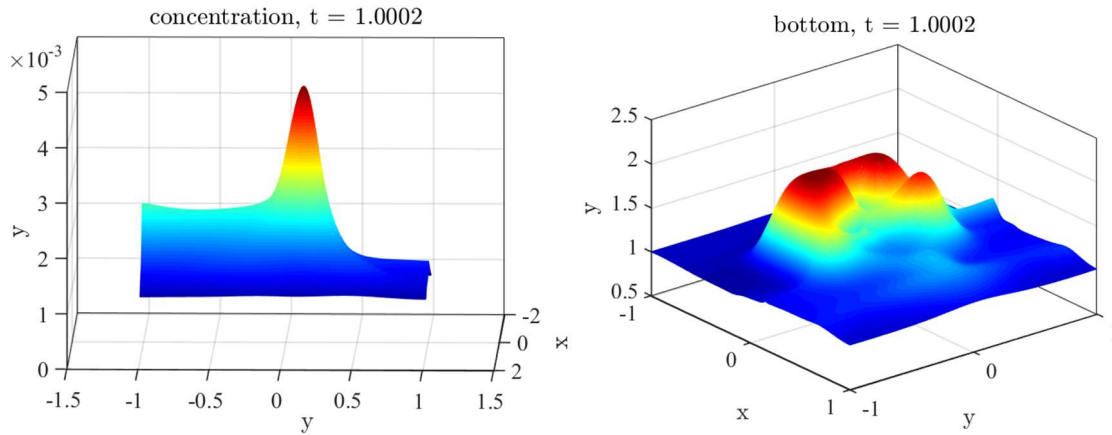

**Figure: 11** Solution of the Riemann problem sediment concentration evolution and morphodynamic profiles, CFL=0.5.

## V.     Conclusion and perspectives

A two-dimensional sediment transport model in nonhomogeneous shallow water equations has been proposed in this work. The model integrates a phase lag effect via a new alternative bedload equation which does not appear in some other models existing in the literature. Moreover, the model captures well the bed wave and the resonance condition via this model can be easily expressed. We proposed an existence theorem of global weak solutions of the model and a convergence study is discussed. It was proved that with this alternative formulation of the bedload equation, the model is still hyperbolic and the finding of the total eigenstructure becomes easy. A new well-balanced positive finite volume method for a 2D sediment transport model has been proposed to solve coastal engineering problems. This method can be applied to several 2D sediment transport models without major modifications. 2D AENO nonlinear reconstruction and second-order Strang method have been presented to obtain second-order accuracy of the fully discrete scheme. High-order accuracy can be simply obtained by increasing the order of derivatives. Considerable attention is paid to the validation of the proposed model by comparing its solutions with experimental data found in the literature. It has proven that our model gives the best result and would need more attention. It is shown that the proposed model describes quite accurately sediment processes even for a large range of sediment diameters. The proposed shock-capturing method can be used for other nonconservative problems associated with other physics without any difficulty. The strategies developed to achieve second order can be modified to other engineering applications or environmental contexts.

**Perspectives**

- A multi-dimensional version of PCCU scheme based on unstructured meshes with mobile domain can be easily design using the same approach.
- Another problem encountered here in the design of the 2D scheme is that the fluxes are computed only at interfaces of the cells and do not account the fluxes at the vertex of each cell. Therefore, it necessary to design a two-dimensional PCCU scheme more general than the current scheme.
- A high order scheme can also be obtained easily but remains an open problem.
- The hydrostatic reconstruction proposed here can be improved to other applications.
- The proposed scheme can be applied to simulate the flooding with sediment deposition in TONGO BASSA basin located in Douala, Cameroon.

**Data availability**

The data that support the findings of this study are available on request from the corresponding author.

**Conflict of Interests**

The author declares that there is no conflict of interest regarding the publication of this paper.

**Acknowledgment**

The author would like to thank an anonymous referee for giving very helpful comments and suggestions that have greatly improved this paper.

# References

[1] Z. Cao, R. Day et S. ,. Egashira, «Coupled and uncoupled numerical modelling of flow and morphological evolution in alluvial rivers,» *Journal of Hydraulic Engineering,* vol. 128, n° %13, pp. 306-321, 2002.

[2] «Bedload transport in shallow water models: Why splitting (may) fail, how hyperbolicity (can) help.,» *Advances in Water Resources,* vol. 34 (8), p. 980{989, 2011.

[3] B. Greimann et J. Huang, «Two-dimensional total sediment load model equations.,» *Haudraul Eng,* vol. 134, pp. 1142-1146, 2008.

[4] L. Xin, A. Mohammadian, Kurganov et J. A. I. Sedano, «Well-balanced central-upwind scheme for a fully coupled shallow water system modeling flows over erodible bed,» *Journal of Computational Physics,* vol. 300, p. 202–218, 2015.

[5] B. Birnir et J. Rowlett, «Mathematical models for erosion and the optimal transportation of sediment,» *arXiv:2012.07736v1 [math.AP],* 2020.

[6] A. R. Ngatcha, B. Nkonga et A. Njifenjou, «Multi-dimensional Positivity-preserving Well-balanced Path-Conservative Central-Upwind scheme on unstructured meshes for a total sediment transport model,» *hal 03668107,* 2022.

[7] M. Castro, J. M. Gallardo et C. Pares, «High-order finite volume schemes based on reconstruction o´f states for solving hyperbolic systems with nonconservative products. applications to shallow-water systems.,» *Mathematics of Computation,* vol. 75, p. 1103–1134, 2006.

[8] N. A. R. Ngatcha, Y. Mimbeu, R. Onguene, S. Nguiya et A. Njifenjou, «A Novel Sediment Transport Model Accounting Phase Lag Effect. A Resonance Condition.,» *WSEAS Transactions on Fluid Mechanics, ISSN / E-ISSN: 1790-5087,* vol. 17, pp. 189-211, 2022.

[9] A. Siviglia, B. Vanzo et E. Toro, «A splitting scheme for the coupled Saint-Venant-Exner model,» *Journal of Adv Res,* vol. 159, p. 104062., 2021.

[10] E. Alqasimi, K. Tew et K. Mahdi, «A new one-dimensional numerical model unsteady hydraulic of sediments in rivers,» *SN APPLIED SCIENCES ,* vol. 2, n° %11480, 2020.

[11] H. P. Gunawan, «Numerical simulation of shallow water equations and related models,» NNT: 2015PEST1010. tel-01216642v2, Paris, 2015.

[12] F. Benkhaldoun, S. Saida et M. Seaid, «A flux-limiter method for dam-break flows over erodible sediment beds,» *Applied Mathematical Modelling,* vol. 36, n° %12012, p. 4847–4861, 2012.

[13] M. J. Castro, E. D. Fernandez-Nieto, T. Morales, G. Narbona-Reina et C. Pares, «A HLLC scheme for nonconservative hyperbolic problems. Application to turbidity currents with sediment transport».

[14] A. Harten, P. Lax et V. Leer, «Upstream differencing and Godunov-type scheme for hyperbolic conservation laws, Upwind and High-Resolution Schemes,» pp. 53-79, 1982.

[15] E. Audusse, C. Chalons et P. Ung, «A simple three-wave Approximate Riemann Solve for the Saint-Venant–Exner equations,» *Communication in mathematics,* 2012.

[16] E. F. Toro, M. Spruce et W. Speares, «Restoration of the contact surface in the HLL-Riemann solver,» *Shock Waves,* vol. 4(1), pp. 25-34, July 1994.

[17] B. Einfeldt, «On Godunov-type methods for dynamics gas,» *SIAM J. Numer. Anal.,* vol. 25, pp. 294-318, 1988.

[18] K. A. Schneider, J. M. Gallardo, S. B. Dinshaw, B. Nkonga et C. Parés, «Multidimensional approximate Riemann solvers for hyperbolic nonconservative systems. Applications to shallow water systems,» *Journal of Computational Physics,* vol. 444, p. 110547, 2021.

[19] M. Castro et E. Fernandez-Nieto, «A class of computationally fast first order finit´e volume solvers: PVM methods,» *SIAM Journal of Scientific Computing,* vol. 34, p. A2173–A2196, 2012.

[20] M. J. Castro, J. M. Gallardo et A. Marquina, «A class of incomplete Riemann solvers based on uniform rational approximations to the absolute value function,» *J. Sci. Comput,* vol. 60, p. 363–389, 2014.

[21] P. C, «Numerical methods for nonconservative hyperbolic systems : a theoretical framework,» *SIAM Journal on Numerical Analysis,* vol. 44, p. 300–321, 2006.

[22] M. Castro Diaz, A. Kurganov et T. Morales de Luna, «PATH-CONSERVATIVE CENTRAL-UPWIND SCHEMES FOR NONCONSERVATIVE HYPERBOLIC SYSTEMS,» *ESAIM: Mathematical Modelling and Numerical Analysis,* vol. 53, pp. 959-985, 2019.

[23] A. Ngatcha et A. Njifenjou, «A well balanced PCCU AENO scheme for a sediment transport model,» *Ocean System Engineering,* vol. 12 https://doi.org/10.12989/ose.2022.12.3, n° %13, 2022.

[24] A. Ngatcha et N. Boniface, « Sediment transport in sheared shallow water flow model: mathematical derivation and hyperbolicity study.,» *In preparation,* 2022.

[25] B. d. A. Saint-Venant, «Théorie du mouvement non-permanent des eaux avec application aux crues des rivières et à l'introduction des marées dans leur lit,» *Comptes Rendus de L'Académie des sciences,* vol. 73, pp. 147-154, 1871.

[26] P. Tassi et C. Villaret, Sisyphe v6.3 User's Manual. EDF R&D., Chatou, France, 2014.

[27] E. F. Toro, A. Santaca, G. I. Montecinos et L. O. Muller, «AENO: a novel reconstruction method in conjunction with ADER schemes for hyperbolic equations.,» *Communications on Applied Mathematics and Computation.,* 2021.

[28] L. Fraccarollo et C. H, «Riemann wave description of erosional dam-break flows,» *J. Fluid Mech. ,* vol. 461, pp. 183-228, 2002.

[29] G. Putu, «Numerical simulation of shallow water equations and related models,» Phd Thesis, Université de Paris-Est, HAL Id: tel 01216642, Paris-Est.

[30] F. Benkhaldoun, S. Saida et M. Seaid, «A flux-limiter method for dam-break flows over erodible sediment beds,» *Applied Mathematical Modelling,* vol. 36, n° %12012, p. 4847–4861, 2012.

[31] V. Mélanie, J. Armelle, M. François et L. B. Sophie, «Experimental Study on Sediment Supply-Limited Bedforms in a Coastal Context,» chez *Sixth International Conference on Estuaries and Coasts (ICEC-2018)*, Caen, France, August 20-23, 2018.

[32] L. C. Rijn, «Sediment Transport, Part II: Suspended Load Transport.,» *Journal of Hydraulic Engineering,* vol. 110, p. 1613—1641, 1984.